\documentclass[letterpaper, 10 pt, conference]{ieeeconf}  
\IEEEoverridecommandlockouts                              

\usepackage{subfig}
\usepackage{tikz,pgf}
\usepackage{graphics} 
\usepackage{epsfig} 
\usepackage{amsmath} 
\usepackage{amssymb}  
\usepackage{color}
\usepackage{relsize}
\usepackage[margin=0.8in]{geometry}

\newcommand{\R}{\mathbb R}

\newtheorem{theorem}{Theorem}

\newtheorem{definition}[theorem]{Definition}

\title{\LARGE \bf
Hysteresis in the Linearized Landau-Lifshitz Equation
}

\author{A. Chow$^{1}$ and K.A. Morris$^{1}$
\thanks{$^{1}$Department of Applied Mathematics, University of Waterloo, Waterloo, Canada, N2L 3G1, {\tt  a29chow@uwaterloo.ca,  kmorris@uwaterloo.ca}. This  research was supported by a Discovery Grant from the Natural Sciences and Engineering Research Council of Canada (NSERC). }}

\begin{document}

\maketitle
\thispagestyle{empty}
\pagestyle{empty}

\begin{abstract}
The Landau-Lifshitz equation describes the behaviour of magnetization inside a ferromagnetic object. It is known that the Landau-Lifshitz equation has an infinite number of stable equilibrium points. The existence of multiple stable equilibria is closely related to hysteresis. This is a phenomenon that is often characterized by a looping behaviour; however, the existence of a loop is not sufficient to identify hysteretic systems, but is defined more precisely as the presence of looping as the frequency of the input goes to zero. We describe these two approaches to identification of  hysteresis and demonstrate that both the linear and nonlinear Landau-Lifshitz equations exhibit hysteresis. The presence of hysteresis in the linear Landau-Lifshitz equation, as well as in a simpler system also described here, indicates that nonlinearity is not necessary for hysteresis to exist.
\end{abstract}

\section{INTRODUCTION}

The Landau-Lifshitz equation describes  magnetic behaviour within ferromagnetic nanostructures such as nanowires \cite{Carbou2006, CarbouLabbe2006, Carbou2008, Gou2011, Noh2012} and nanoparticles \cite{Mayergoyz2010}.  These structures are  used in memory storage devices such as hard disks or tape recordings. 
Systems modelled by the Landau-Lifshitz equation are generally regarded as exhibiting hysteresis \cite{CarbouEfendiev2009, Wiele2006, Wiele2007, Yang2011}. Hysteresis is a phenomenon that occurs in many processes \cite{Morris2011}. Examples include other types of  magnetization \cite{Bertotti1998,  Cowburn1999, Morris2011, Schneider2006}, smart materials  \cite{Smith2005}, \cite{Valadkhan2009} freezing and thawing processes  \cite{Alimov2002, Alimov1998, Christenson2001,Petrov2006}, and predator-prey relationships \cite{Aiki2005},  \cite[Section~1.3]{Murray1993}.  

It is difficult to define hysteresis precisely.
In much of the literature, hysteresis is represented by hysteresis operators \cite{Aiki2005, Brokate1996, CarbouEfendiev2009, Logemann2008, Valadkhan2010, Vistinin1983, Visintin2006}. The definition of hysteresis operator can include operators that are not commonly understood as hysteretic \cite{Morris2011}.  Furthermore, the notion of a hysteretic operator is not appropriate for the Landau-Lifshitz equation since it requires input-output behaviour to be rate-independent; that is, the input-output diagram depends only on the values of the input and not on its frequency. 

A common theme in defining hysteresis is that of a looping behaviour displayed in the input-output map.   The shape of the loop can appear quite different under different inputs, even for the same system. The Landau-Lifshitz equation is an example of this, which is demonstrated later in the paper. There are other definitions of hysteresis available in the literature, aside from a hysteresis operator.  One such definition defines hysteretic systems as having multiple stable equilibrium points and dynamics that are faster than the rate at which inputs are varied \cite{Morris2011}.  Another definition requires persistent looping behaviour in the input-output maps  \cite{Bernstein2005}. We describe both approaches and  show that both the linear and nonlinear Landau-Lifshitz equation exhibits hysteresis. 




\section{Hysteresis}\label{secHys}
Hysteresis can be considered from a dynamical systems perspective.
 \begin{definition}\cite[Definition~3]{Morris2011} \label{defmultiequilibrium}
A system exhibits {\it hysteresis} if it has \\
(a) multiple stable equilibrium points and\\
(b) dynamics that are considerably faster than the time scale at which inputs are varied.
\end{definition}
The standard definition of stable equilibrium point \cite[Definition~3.2]{Walker1976} with zero input is used here. Condition~(b) is relative to the speed at which a controlled input is changed. In many cases, hysteresis is present but is rate-dependent \cite{Morris2011}. 

\begin{figure}\vspace{-0.3cm}
  \centering \hspace*{-0.25cm}
  \subfloat[]{\label{figa} \includegraphics[width=0.17\textwidth, trim=0 200 0 200]{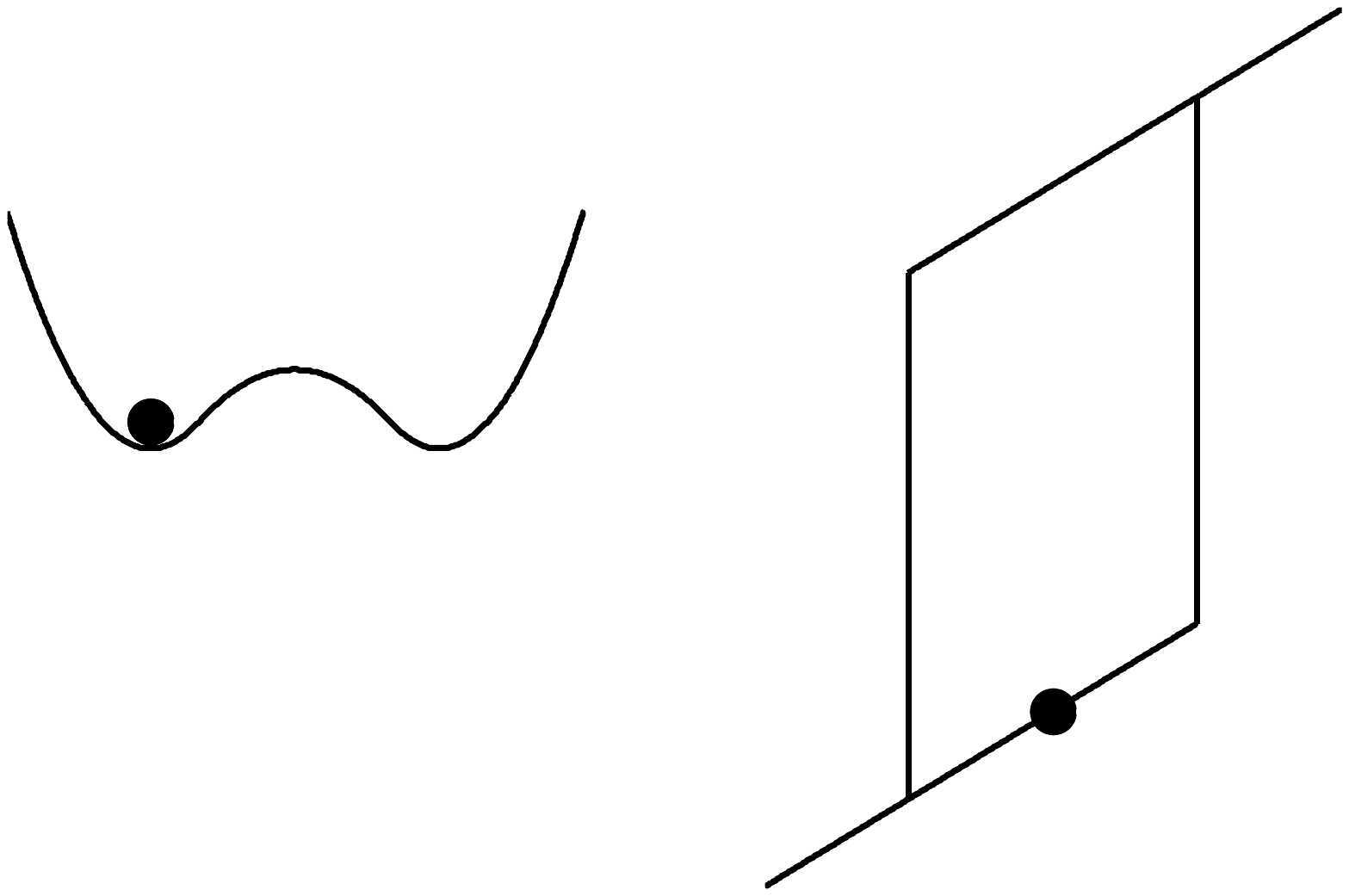}} 
  \subfloat[]{\label{figb}\includegraphics[width=0.17\textwidth, trim=0 200 0 200]{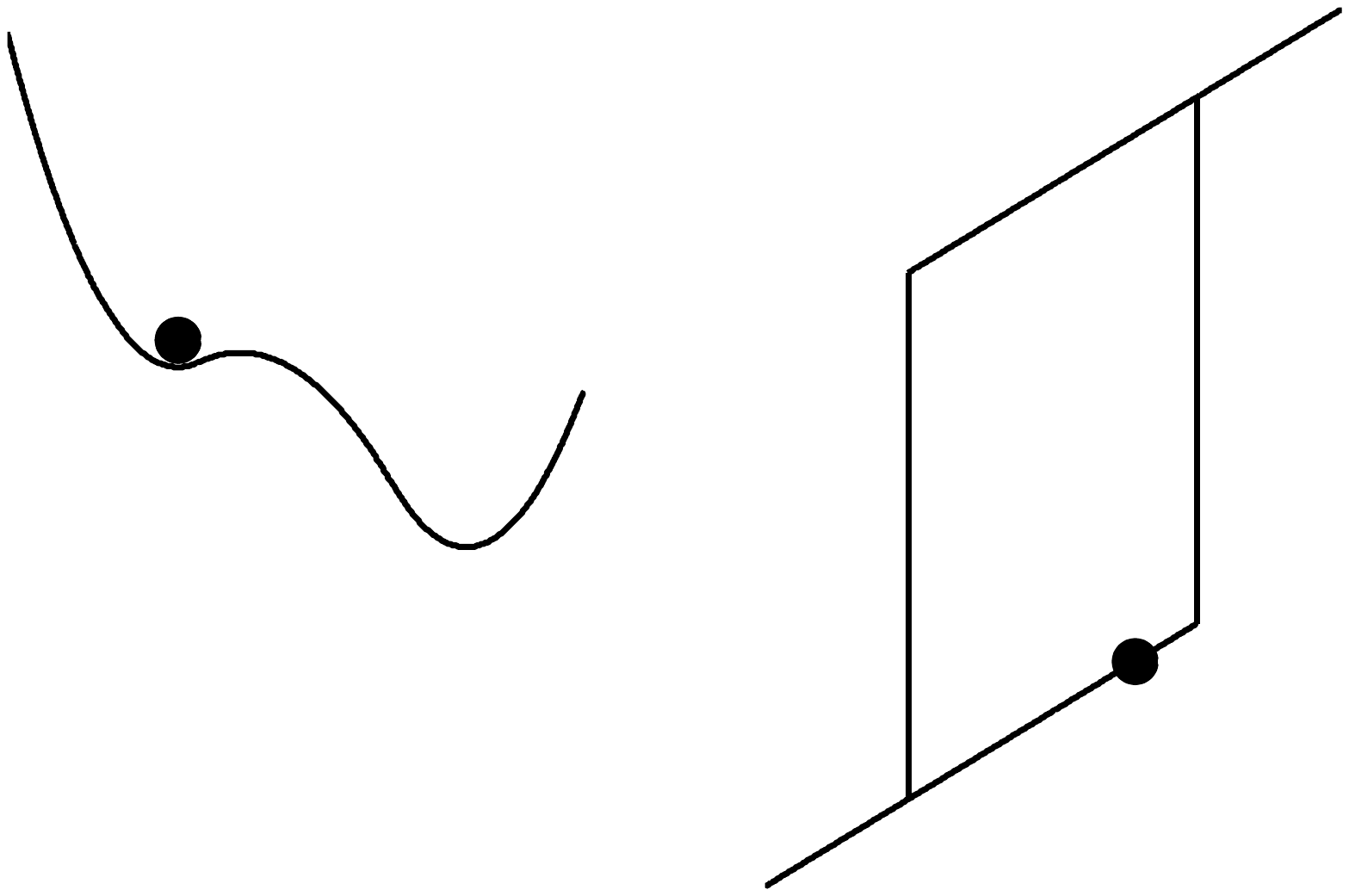}}
    \subfloat[]{\label{figc} \includegraphics[width=0.17\textwidth, trim=0 200 0 200]{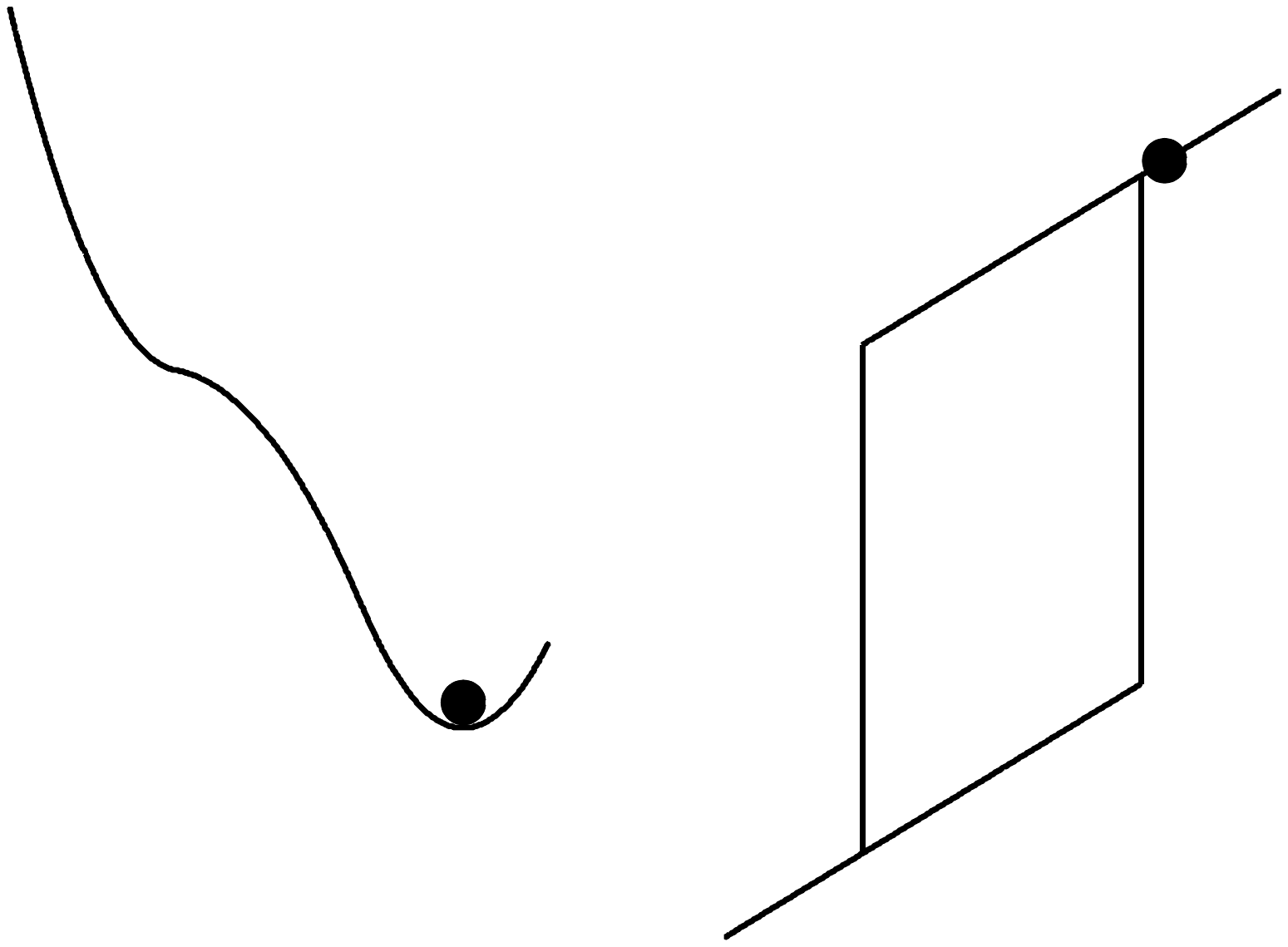}} \\
     \centering  \hspace*{-0.25cm}
  \subfloat[]{\label{figd}\includegraphics[width=0.17\textwidth, trim=0 200 0 200]{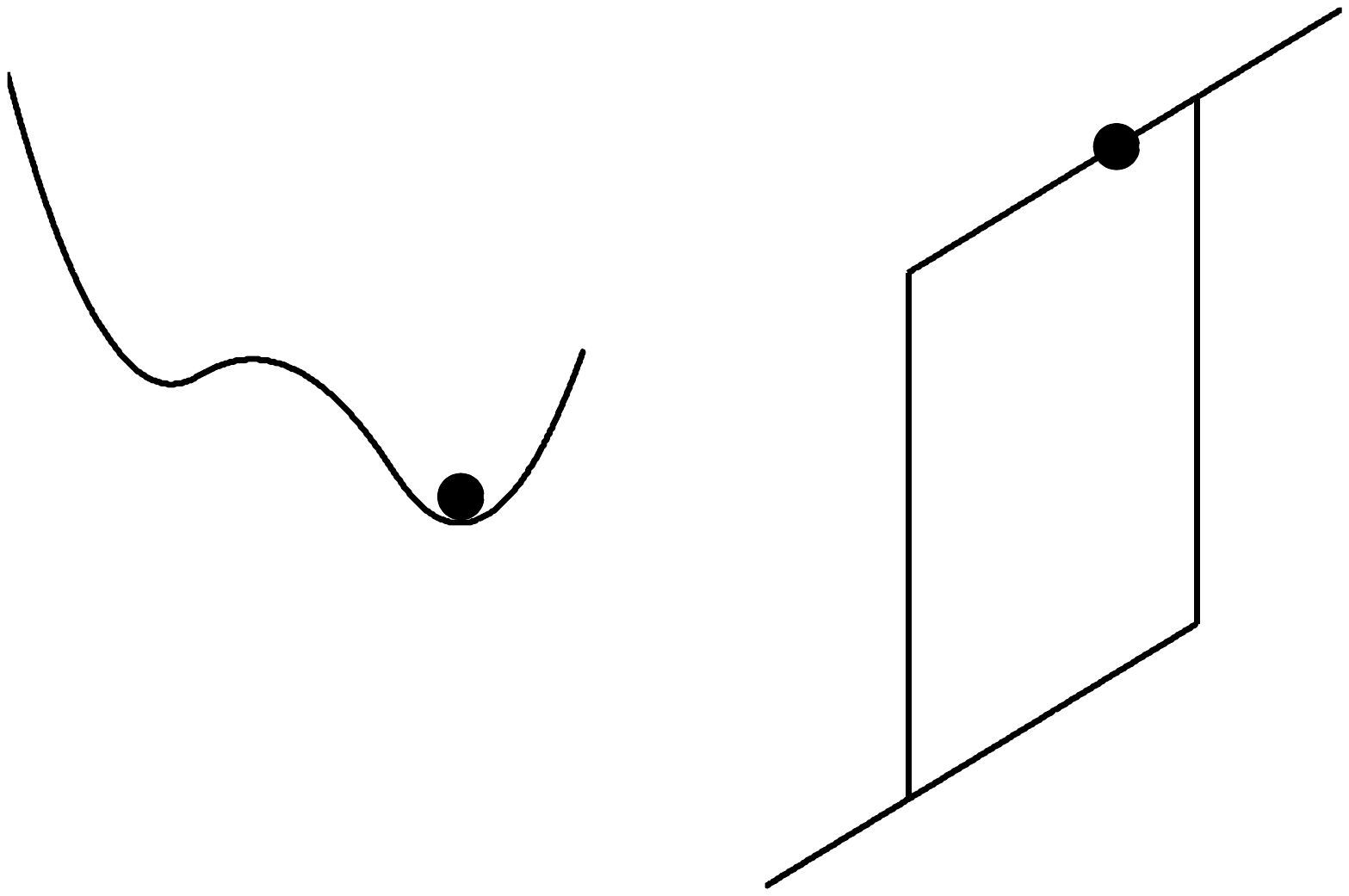}}
  \subfloat[]{\label{fige}\includegraphics[width=0.17\textwidth, trim=0 200 0 200]{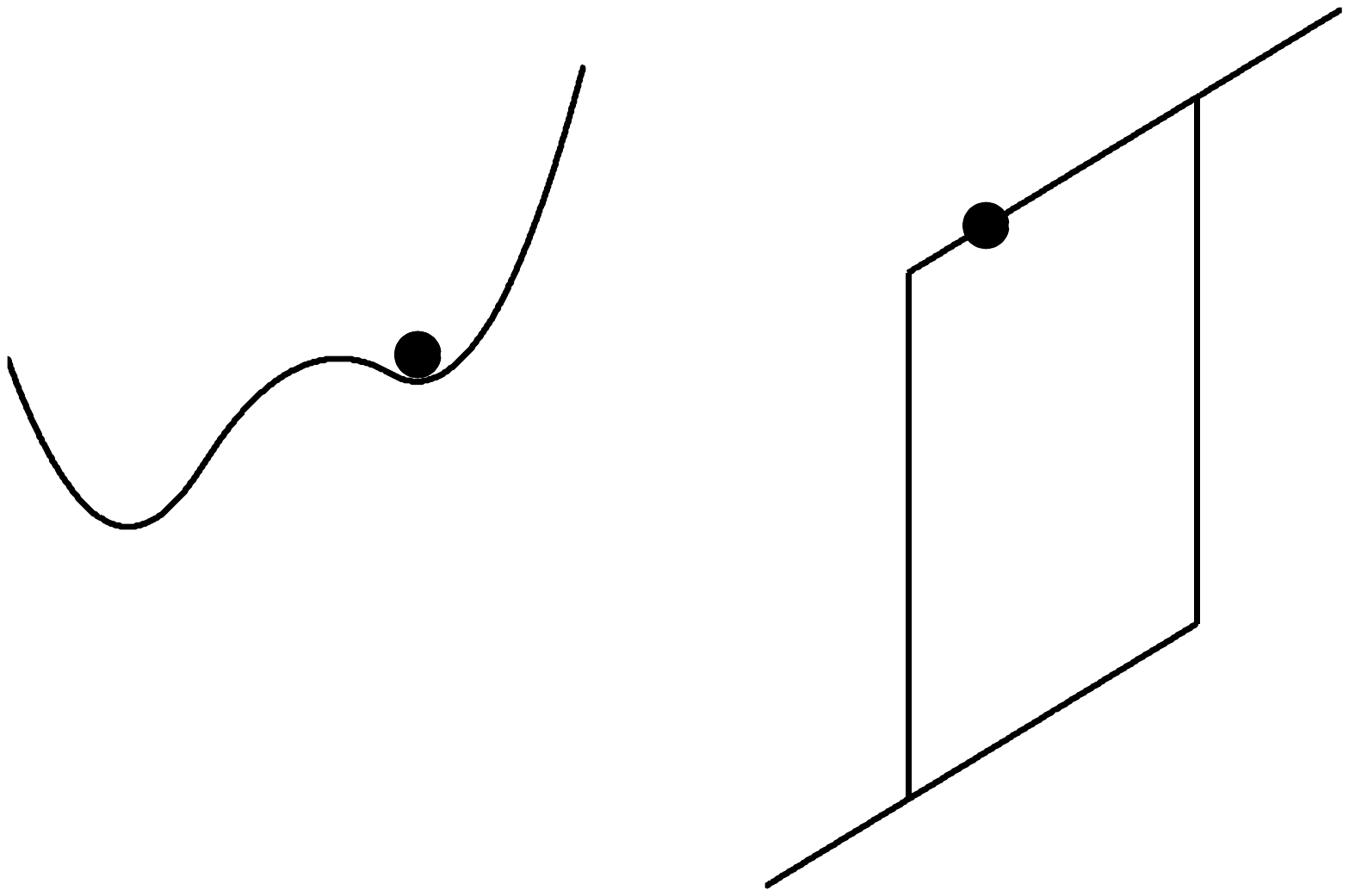}}
  \subfloat[]{\label{figf}\includegraphics[width=0.17\textwidth, trim=0 200 0 200]{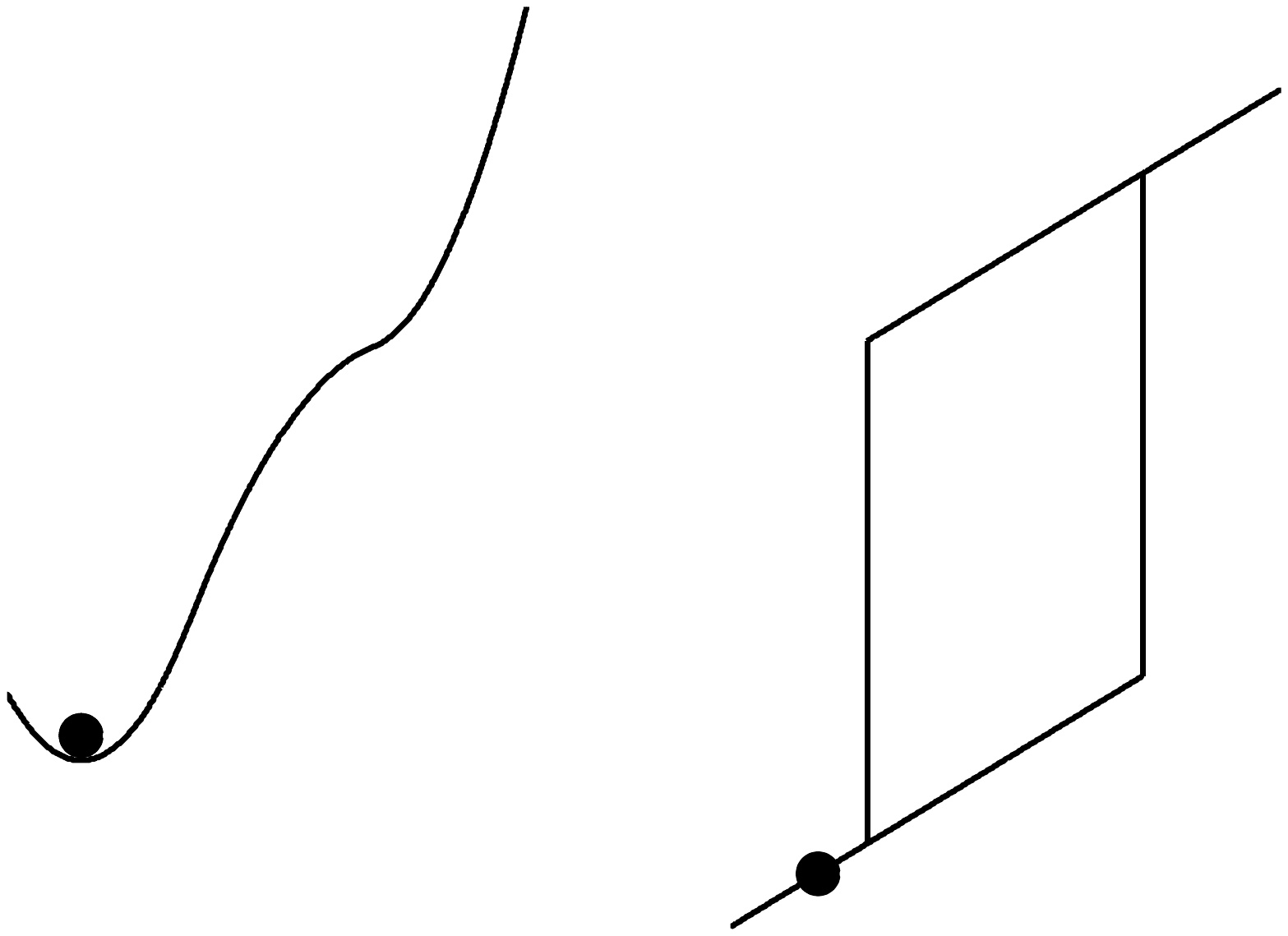}}\\
   \centering  \hspace*{-0.25cm}
  \subfloat[]{\label{figg}\includegraphics[width=0.17\textwidth, trim=0 200 0 200]{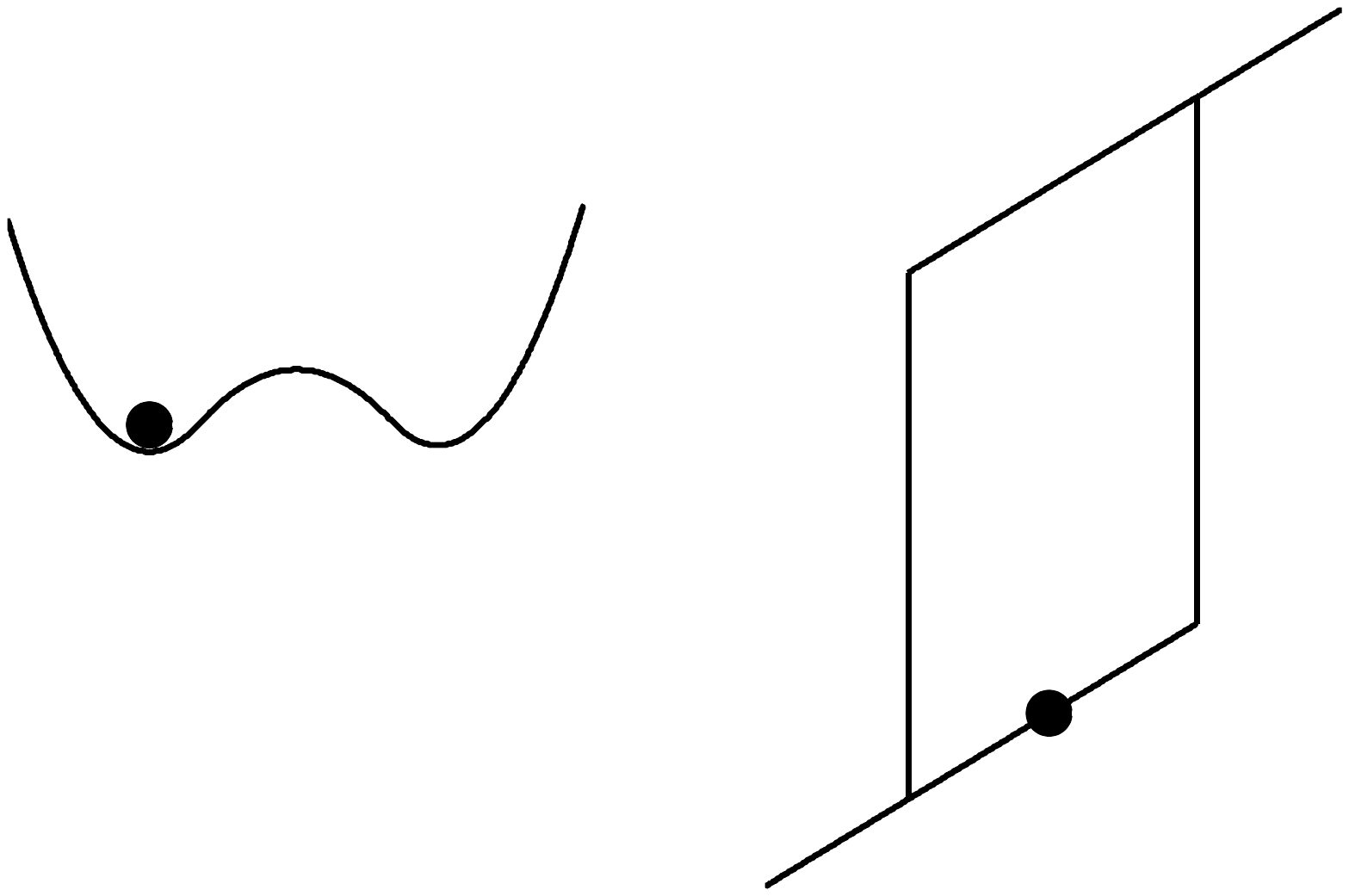}} 
    \subfloat{\includegraphics[width=0.17\textwidth, trim=0 200 0 200]{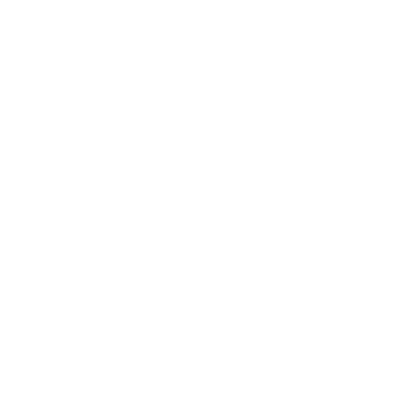}}
        \subfloat{\includegraphics[width=0.17\textwidth, trim=0 200 0 200]{blank.png}}
  \caption{\label{figmultiequilibrium}  \small  Dynamics of a system with two stable equilibria for a range of values of the input, and  corresponding hysteresis loop.  Subfigures (a)-(g) show the system variables as the input changes with time.   The left-hand curve in each subfigure is the energy functional with the horizontal axis as the system state and the vertical axis as the energy. The right-hand curve in each subfigure is the system input-output diagram with input  on the horizontal axis and  output on the vertical axis.}
\end{figure}

Definition~\ref{defmultiequilibrium} is related to the looping behaviour often associated with hysteresis.  Consider  a system with two stable equilibria and suppose the system is initially at the left equilibrium (Figure~\ref{figa}).  If the input is increased, the system will tend to stay in equilibrium (Figure~\ref{figb}) with only a small move upward along the hysteresis curve.  When the input is increased enough such that the equilibrium disappears, the system moves to the right equilibrium (Figure~\ref{figc}).    This corresponds to moving along the steepest portion of the hysteresis loop.  The system stays at the right equilibrium (Figure~\ref{figd}). If the input decreases enough so that the right equilibrium disappears (Figure~\ref{fige}), the system moves back to the left equilibrium (Figures~\ref{figf}, \ref{figg}).  
 For systems which move to equilibrium faster than changes in the input, the transition from one equilibrium to the other is nearly instantaneous and the system behaviour can appear rate-independent. For further discussion, see \cite{Bernstein2007,Morris2011}.

However, the existence of a loop in the input-output map is not sufficient to define a system as hysteretic. Consider the simple example of a damped second-order system
\begin{equation}\label{eqspringsystem}
\ddot{y}(t)+c\dot y(t)+ky(t)=u(t).
\end{equation}
 Writing (\ref{eqspringsystem})  in first-order form with $u(t)\equiv0$  and $[x_1,x_2]=[y,\dot y]$ leads to
\begin{subequations}\label{eqspringsystemfirstorder}
\begin{align}
\dot x_1(t)&=x_2(t)\\
\dot x_2(t)&{=-cx_2(t)-kx_1(t)}
\end{align}
\end{subequations}
Analyzing the unforced system, the  only equilibrium point is $(0,0)$ and the eigenvalues are
\begin{align}\label{eqeigenvalues}
\lambda_{1,2} &=\frac{-c\pm\sqrt{c^2-4k}}{2}.
\end{align}
If both eigenvalues have nonpositive real part, then the equilibrium point $(0,0)$ is stable.  

For arbitrary initial conditions, $y(0)=y_0$ and  $\dot{y}(0)=y_1$, the solution to (\ref{eqspringsystem}) with $u(t)=\sin(\omega t)$ 
 is
\begin{align*}
y(t)&=\frac{1}{\sqrt{221}}\left( \frac{\left(-\lambda_1-\lambda_2\omega^2\right)y_0+\omega}{\omega^2-1-15\lambda_1}+y_1\right) e^{\lambda_1t}\\
&-\frac{1}{\sqrt{221}}\left( \frac{\left(-\lambda_2-\lambda_1\omega^2\right)y_0+\omega}{\omega^2-1-15\lambda_2}+y_1\right) e^{\lambda_2t}\\
&+\frac{15\omega\cos(\omega t)+(\omega^2-1)\sin(\omega t)}{2\omega^2-1-\omega^2(15^2+\omega^2)}
\end{align*}
where $\lambda_{1},\lambda_{2}$ are defined in (\ref{eqeigenvalues}).
It follows that 
\begin{align*}
\lim_{\omega\rightarrow 0} y(t)&=\frac{1}{\sqrt{221}}\left( \frac{\lambda_1y_0}{1+15\lambda_1}+y_1\right) e^{\lambda_1t}\\
&-\frac{1}{\sqrt{221}}\left( \frac{\lambda_2y_0}{1+15\lambda_2}+y_1\right) e^{\lambda_2t}.
\end{align*}
If both the real parts of $\lambda_{1}$ and $\lambda_{2}$ are negative, then $y(t)\rightarrow0$ as $\omega \rightarrow 0$, $t\rightarrow \infty$, which means regardless of the initial condition, the steady state is zero as $\omega \rightarrow 0$ and $t\rightarrow \infty$.
  
\begin{figure}\vspace*{-0.5cm}
  \centering \hspace*{-0.25cm}
  \subfloat[$\omega=1$]{\label{} \includegraphics[trim = 0mm 60mm 0mm 60mm, clip, width=0.26\textwidth]{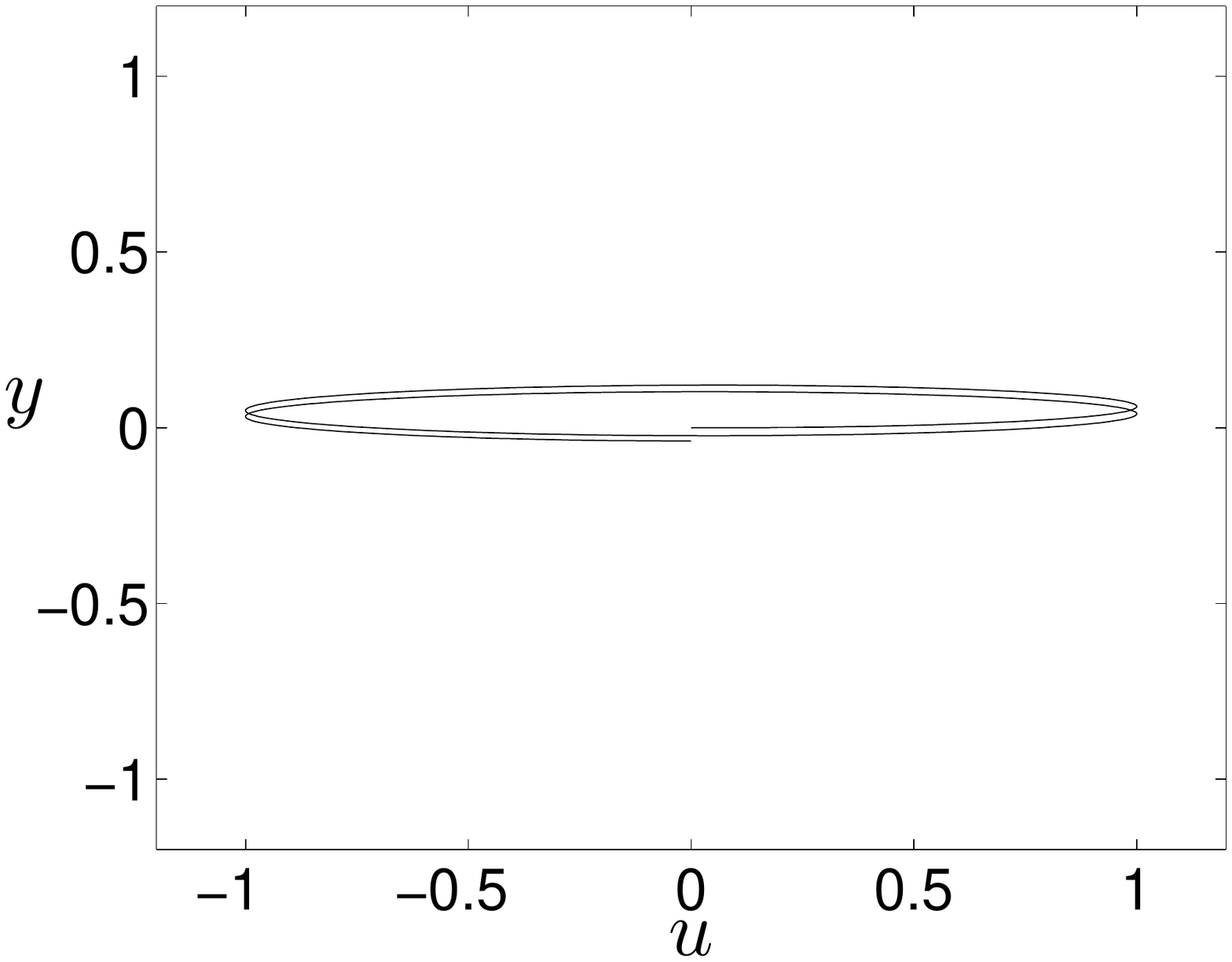}} 
  \subfloat[$\omega=0.1$]{\label{}\includegraphics[trim = 0mm 60mm 0mm 60mm, clip,width=0.26\textwidth]{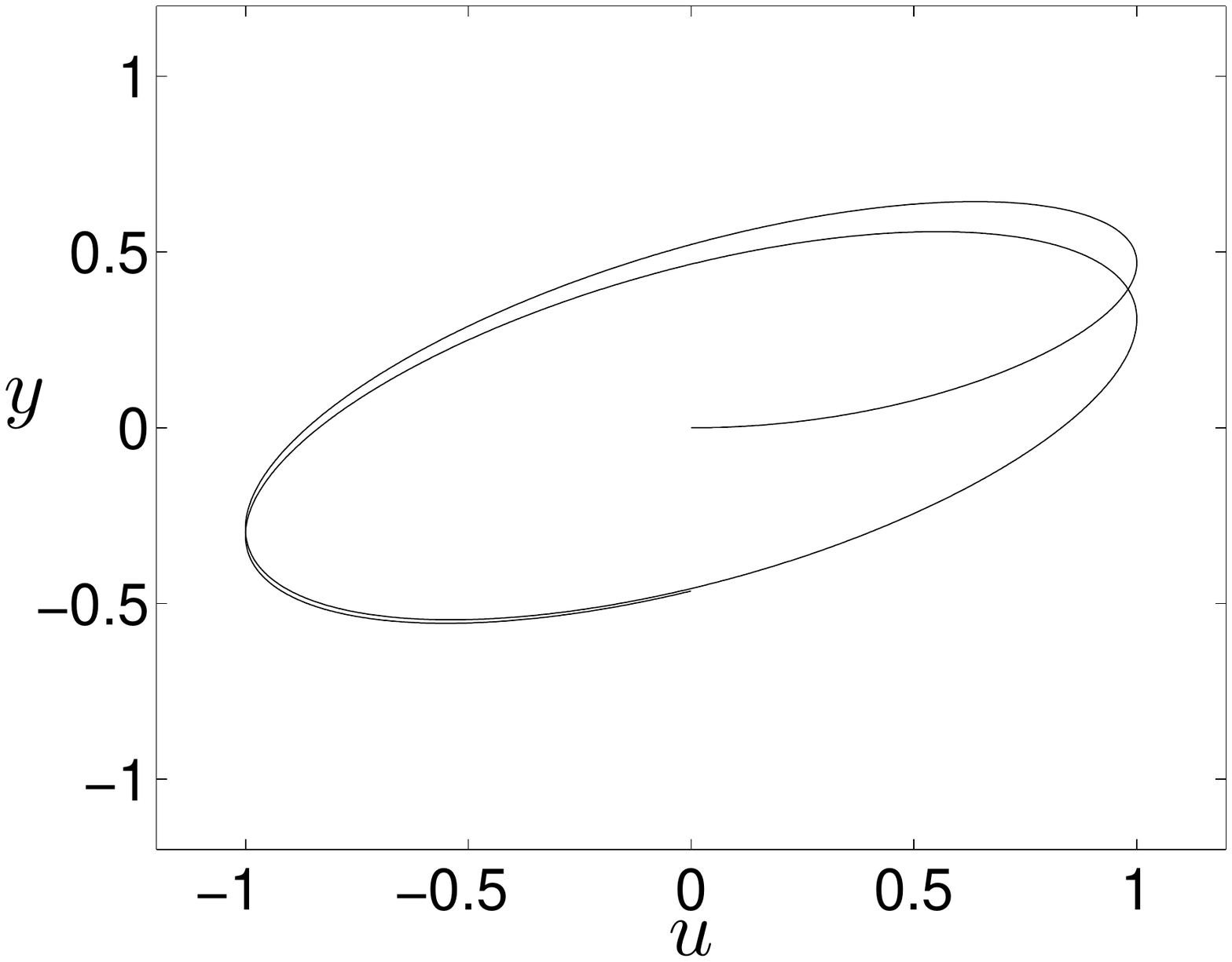}}\\
 \centering  \hspace*{-0.25cm}
    \subfloat[$\omega=0.01$]{\label{} \includegraphics[trim = 0mm 60mm 0mm 60mm, clip,width=0.26\textwidth]{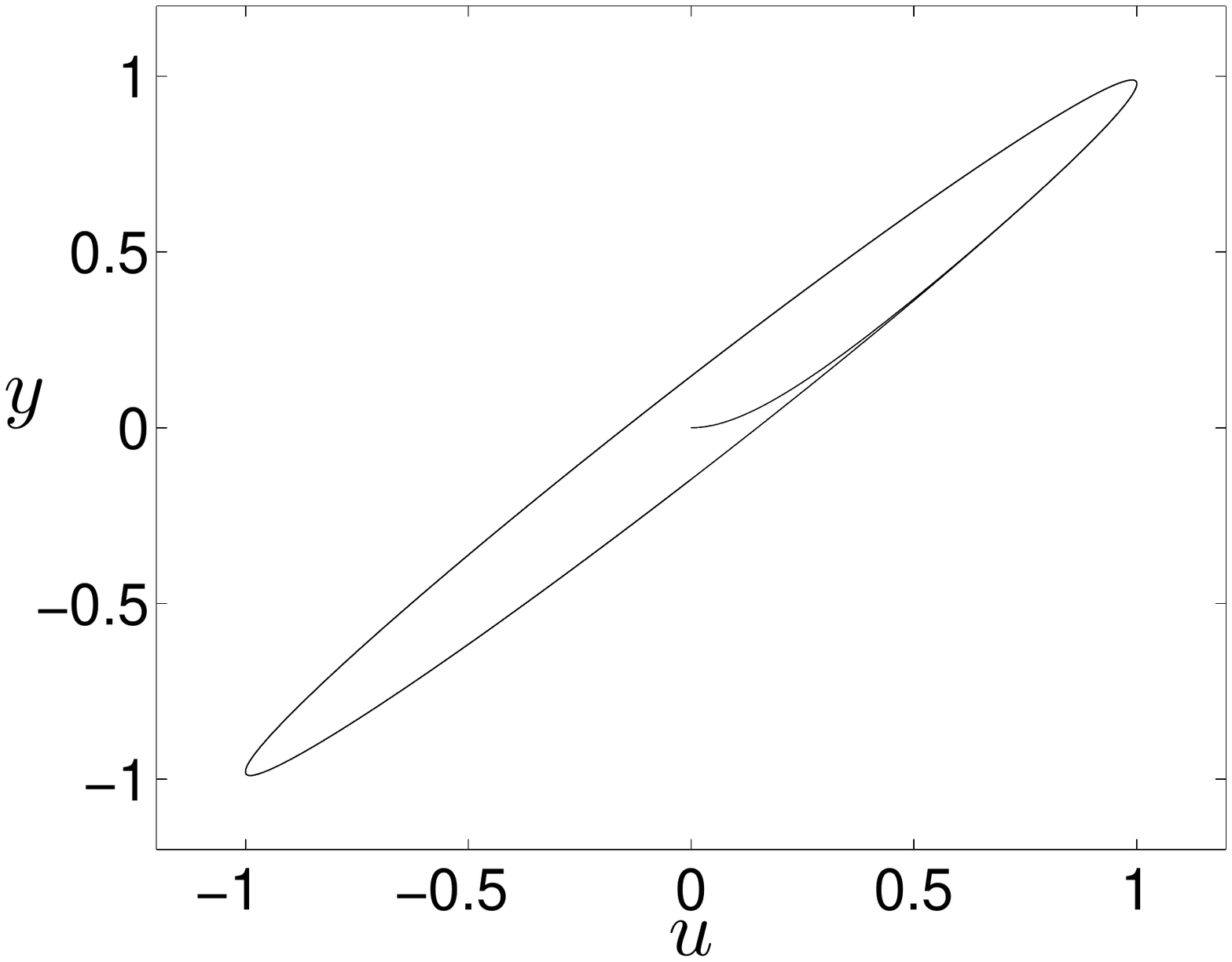}} 
  \subfloat[$\omega=0.001$]{\label{}\includegraphics[trim = 0mm 60mm 0mm 60mm, clip,width=0.26\textwidth]{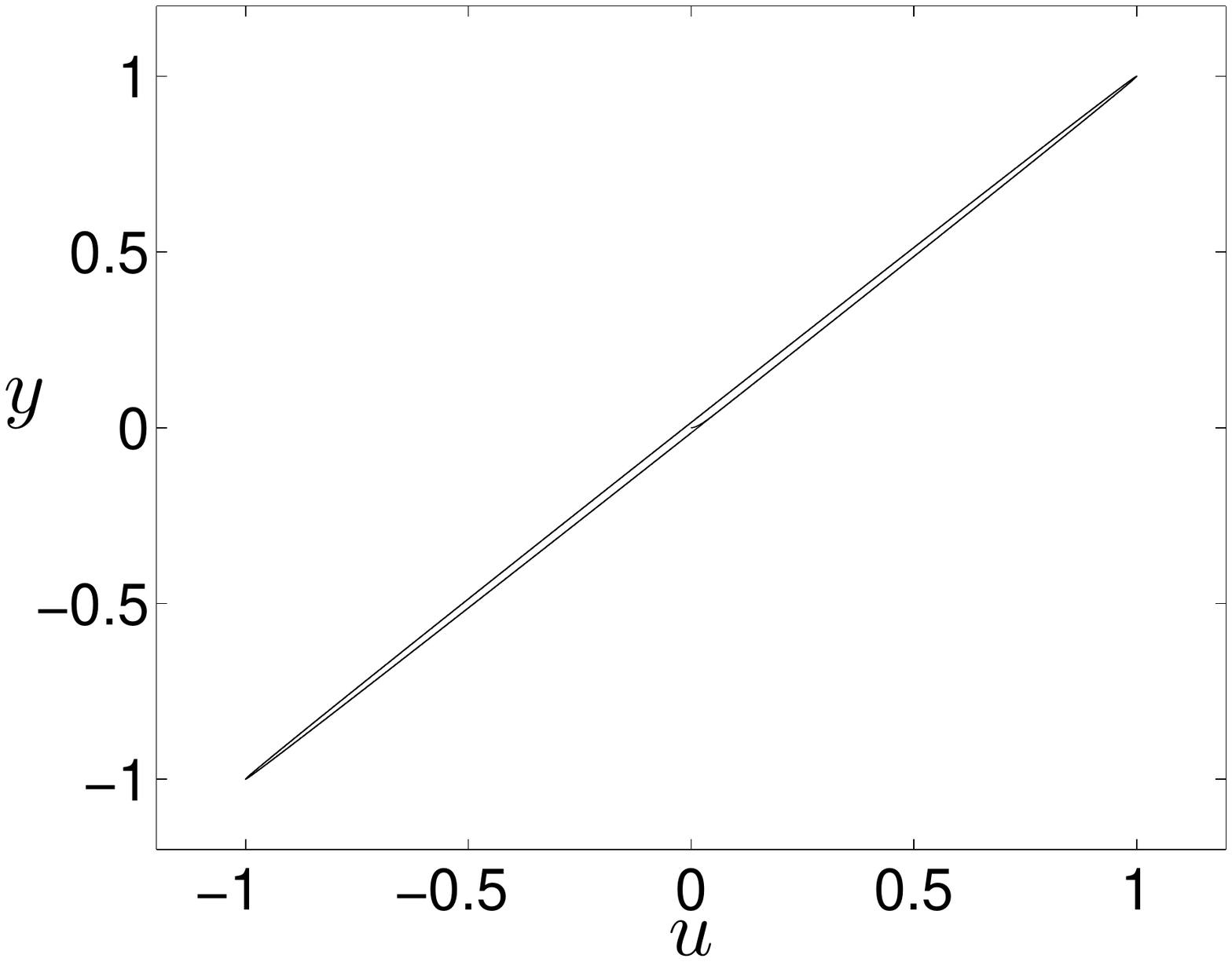}}
  \caption{\label{figspringtrolleyloop} \small Input-output curves for equation~(\ref{eqspringsystem}) with $c=15$ and $k=1$.  The  initial position is $(y,\dot y)=(0,0)$ and the input is $u(t)=\sin(\omega t)$ for various frequencies $\omega$. The uncontrolled system has one equilibrium point. As $\omega \to 0, $ the loops disappear. This system does not display hysteresis.}
\end{figure}

Since~(\ref{eqspringsystemfirstorder}) has exactly one equilibrium, the system does not exhibit hysteresis according to Definition~\ref{defmultiequilibrium}.
Figure~\ref{figspringtrolleyloop} depicts the input-output curves of (\ref{eqspringsystem}), with $c=15$, $k=1$, initial condition $(y(0),\dot y(0))=(0,0)$ and with input $u(t)= \sin(\omega t)$ for different frequencies $\omega .$
   For large $\omega$ there is a loop in the input-output diagram. However,  as $\omega$ approaches zero, the loop  vanishes. This illustrates  a second definition of hysteresis.  
\begin{definition}\label{defBernstein}\cite{Bernstein2005}
A system exhibits {\it hysteresis} if  a nontrivial closed curve in the input-output map persists for a periodic input as the frequency component of the input signal approaches zero. 
\end{definition}
For more details on Definition~\ref{defBernstein}, see also \cite{Ikhouane2013}. Note that a closed curve is said to be non-trivial if it has a non-empty interior. For example, a constant function is a trivial closed curve.

Definitions~\ref{defmultiequilibrium} and~\ref{defBernstein} are complementary, not contradictory.  Definition \ref{defmultiequilibrium} refers to the internal behaviour of the system state while Definition~\ref{defBernstein} refers to the system input-output behaviour. 
Note that equation~(\ref{eqspringsystem}) satisfies neither definition. 

\begin{figure}[h]\vspace*{-0.5cm}
  \centering \hspace*{-0.5cm}
  \subfloat[$\omega=1$]{\label{} \includegraphics[trim = 0mm 60mm 0mm 60mm, clip, width=0.26\textwidth]{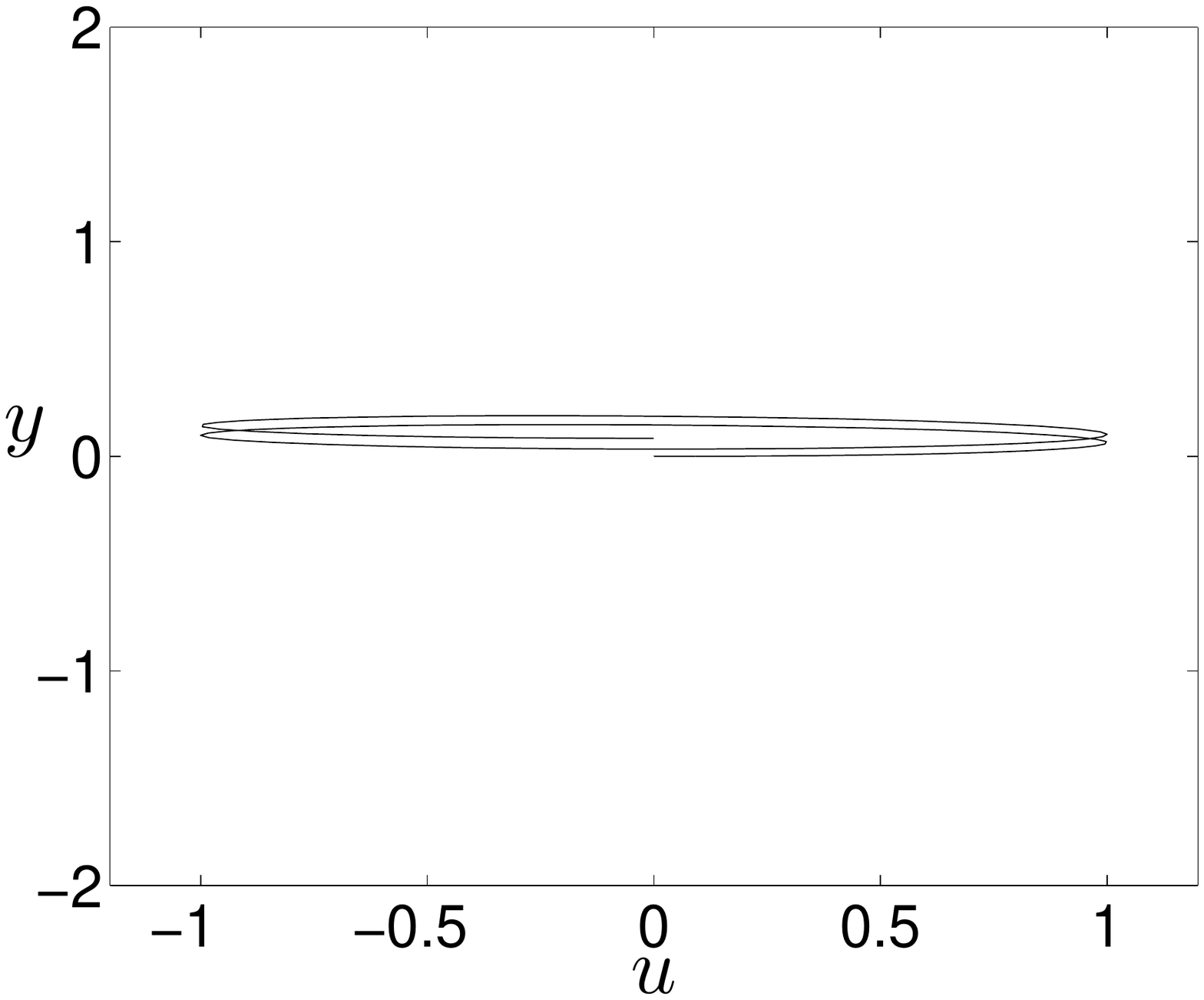}} 
  \subfloat[$\omega=0.1$]{\label{}\includegraphics[trim = 0mm 60mm 0mm 60mm, clip,width=0.26\textwidth]{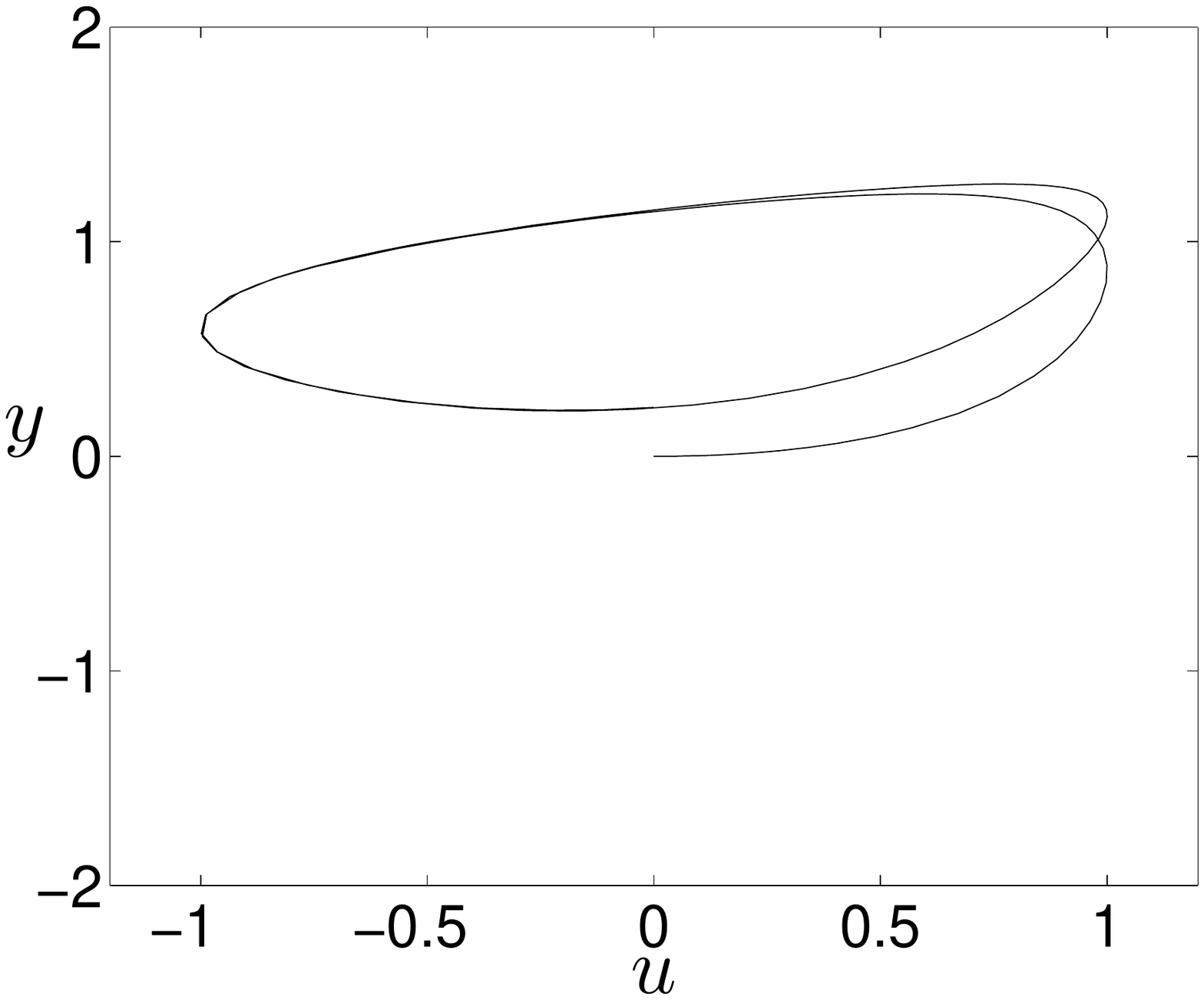}}\\
 \centering  \hspace*{-0.5cm}
    \subfloat[$\omega=0.01$]{\label{} \includegraphics[trim = 0mm 60mm 0mm 60mm, clip,width=0.26\textwidth]{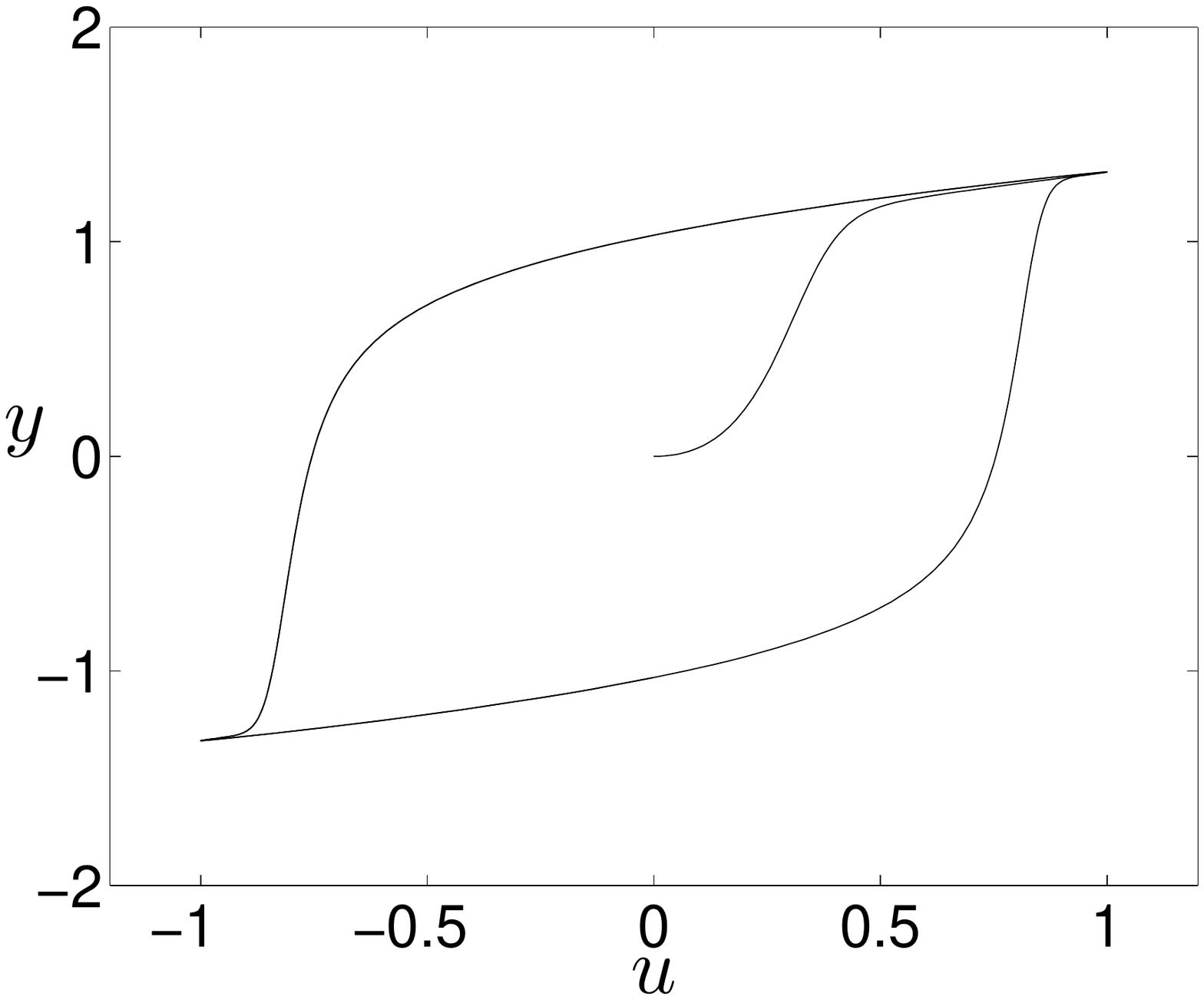}} 
  \subfloat[$\omega=0.001$]{\label{}\includegraphics[trim = 0mm 60mm 0mm 60mm, clip,width=0.26\textwidth]{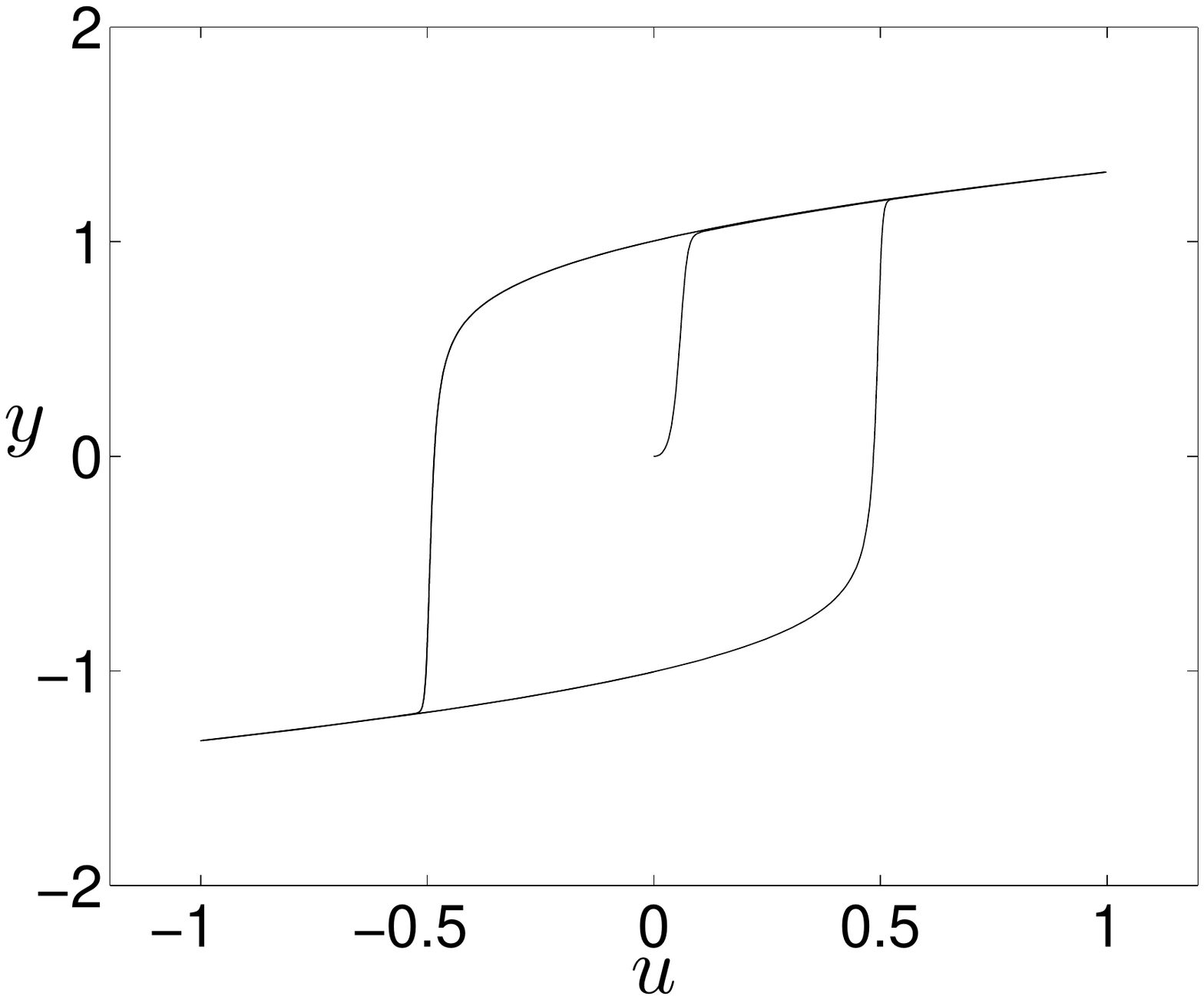}}\\
  \caption{\label{fignonlinearspringtrolleyloop}\small  Input-output curves for the nonlinear system in~(\ref{eqnonlinearspringsystem}) with $c=15$ and $k=-1$. The initial position is $(y,\dot y)=(0,0)$ and the input is $u(t)=\sin(\omega t)$ for various $\omega$.}
\end{figure}

Modify model (\ref{eqspringsystem})  to obtain the nonlinear system 
\begin{equation}\label{eqnonlinearspringsystem}
\ddot{y}(t)+c\dot y(t)+k\left(y(t)-y^3(t)\right)=u(t).
\end{equation}
 Setting $c=15$ and $k=-1$,  the input-output curves with zero initial condition and periodic input are shown in Figure~\ref{fignonlinearspringtrolleyloop}.  It is clear from the figure that even for small $\omega$ a loop persists, which by Definition~\ref{defBernstein} indicates the presence of  hysteresis.  Observe also the differently shaped hysteresis loops for different inputs, which illustrate the complex behaviour of hysteresis, even for this second-order system.

\begin{figure}[h]
\centering\includegraphics[trim = 0mm 60mm 0mm 60mm, clip, width=0.335\textwidth]{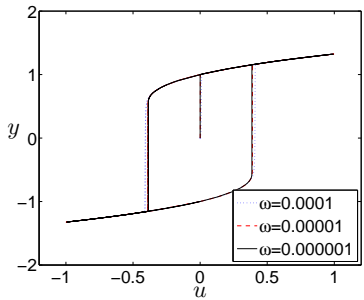}
\caption{\label{figrateindependencespringtrolley} \small Input-output diagrams for the nonlinear system in (\ref{eqnonlinearspringsystem}) with $c=15$, $k=-1$. The initial position is  $(y,\dot y)=(0,0)$ and the input is $u(t)=\sin(\omega t)$ for $\omega\leq 0.0001$. At these low frequencies, the system exhibits rate independence. }
\end{figure}

System  (\ref{eqnonlinearspringsystem}) also satisfies Definition~\ref{defmultiequilibrium}. Rewriting (\ref{eqnonlinearspringsystem}) with no input $u$,
\begin{align*}
\dot x_1(t)&=x_2(t)\\
\dot x_2(t)&=-cx_2(t)-kx_1(t)+kx_1^3(t).
\end{align*}
The equilibrium points are $(0,0)$ and $(\pm 1,0)$. The  eigenvalues of  the system linearized about $(\pm 1,0)$ have negative real part and hence $(\pm 1,0)$ are stable equilibrium points. This satisfies the first condition in Definition~\ref{defmultiequilibrium}. (The equilibrium $(0,0)$ is unstable since one of  its corresponding eigenvalues has positive real part.) Figure~\ref{figrateindependencespringtrolley} demonstrates that the dynamics in equation~(\ref{eqnonlinearspringsystem}) are independent of the rate at which inputs are varied given the range of frequencies.  This satisfies the second condition in Definition~\ref{defmultiequilibrium}. 
A similar example for the displacement of a magnetic beam is presented in  \cite[Example~4]{Morris2011}.  

Consider again  the linear equation~(\ref{eqspringsystem}) but with $k=0$; that is
\begin{equation}\label{eqspringsystemkis0}
\ddot{y}(t)+c\dot y(t)=u(t).
\end{equation} 
and its corresponding unforced system in first-order form is 
\begin{subequations}\label{eqspringsystemfirstorderMultieq}
\begin{align}
\dot x_1(t)&=x_2(t)\\
\dot x_2(t)&={-cx_2(t)}.
\end{align}
\end{subequations}
The unforced system (\ref{eqspringsystemfirstorderMultieq}) has an infinite number of equilibria, $(a,0)$ where $a$ is any constant. The eigenvalues of  (\ref{eqspringsystemfirstorderMultieq}) are $0$ and $-c$ and so the equilibrium points are stable  if $c>0$. 

\begin{figure}[h]\vspace*{-0.3cm}
  \centering \hspace*{-0.5cm}
  \subfloat[$\omega=1$]{\label{} \includegraphics[trim = 0mm 60mm 0mm 60mm, clip, width=0.26\textwidth]{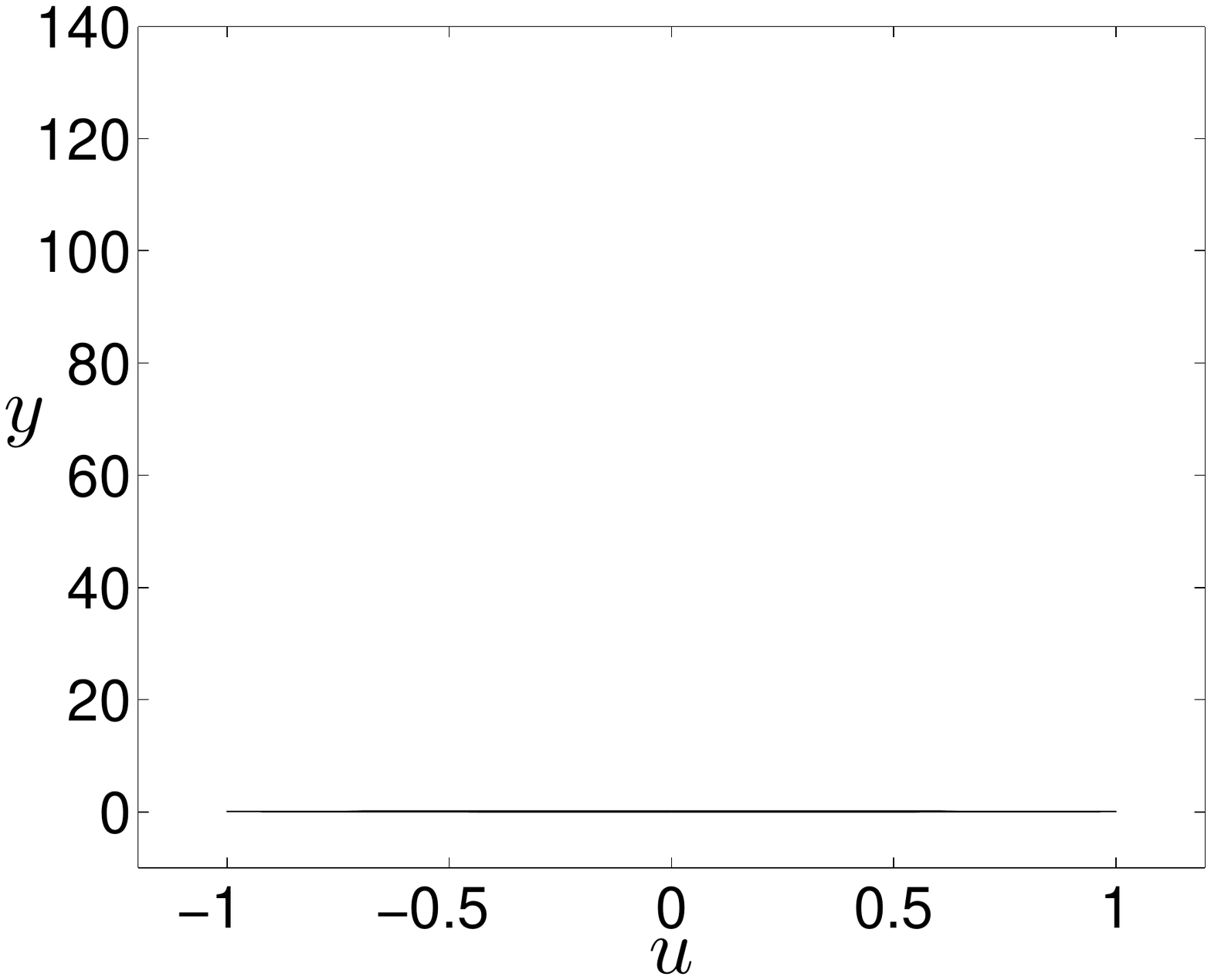}} 
  \subfloat[$\omega=0.1$]{\label{}\includegraphics[trim = 0mm 60mm 0mm 60mm, clip,width=0.26\textwidth]{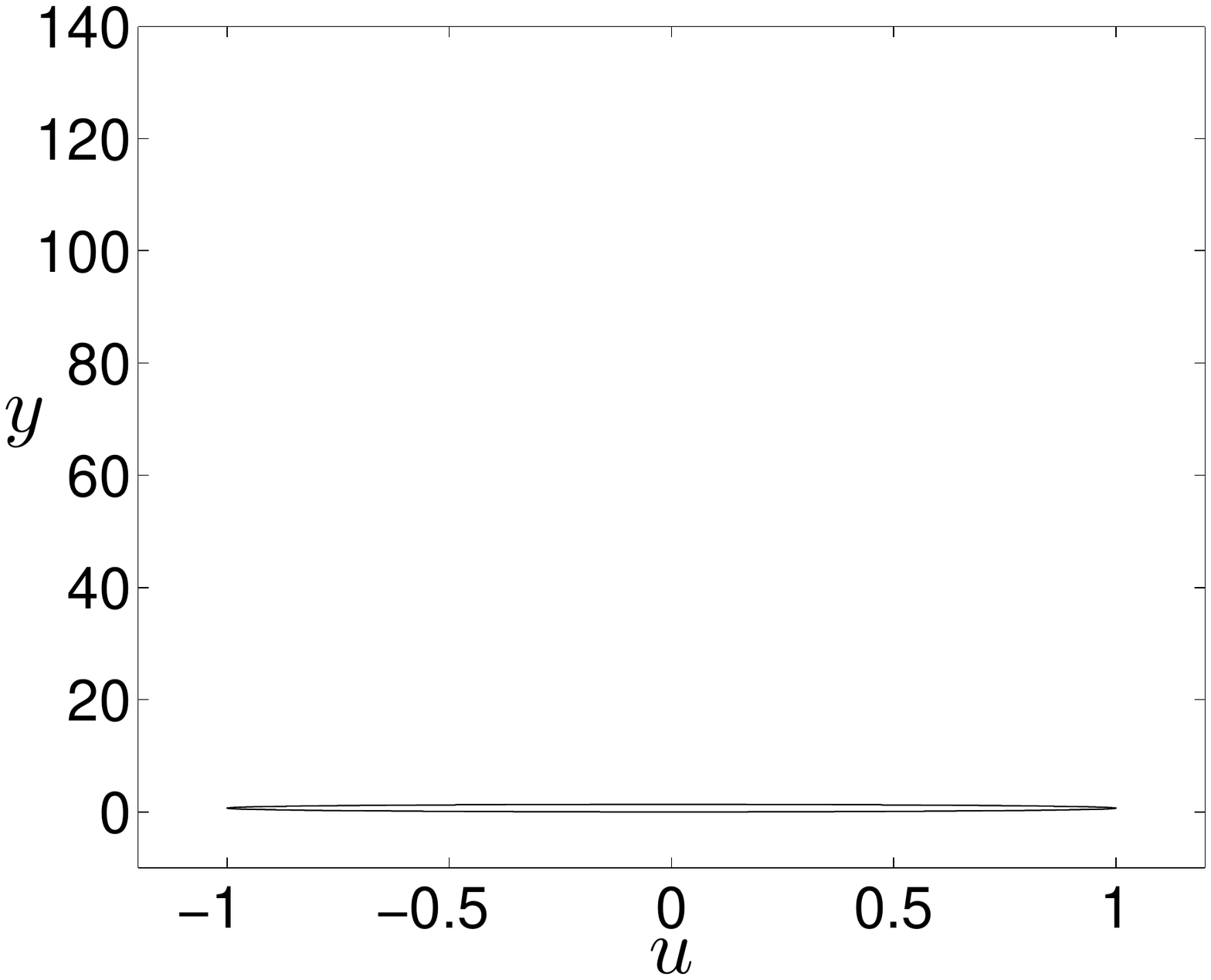}}\\
 \centering   \hspace*{-0.5cm}
    \subfloat[$\omega=0.01$]{\label{} \includegraphics[trim = 0mm 60mm 0mm 60mm, clip,width=0.26\textwidth]{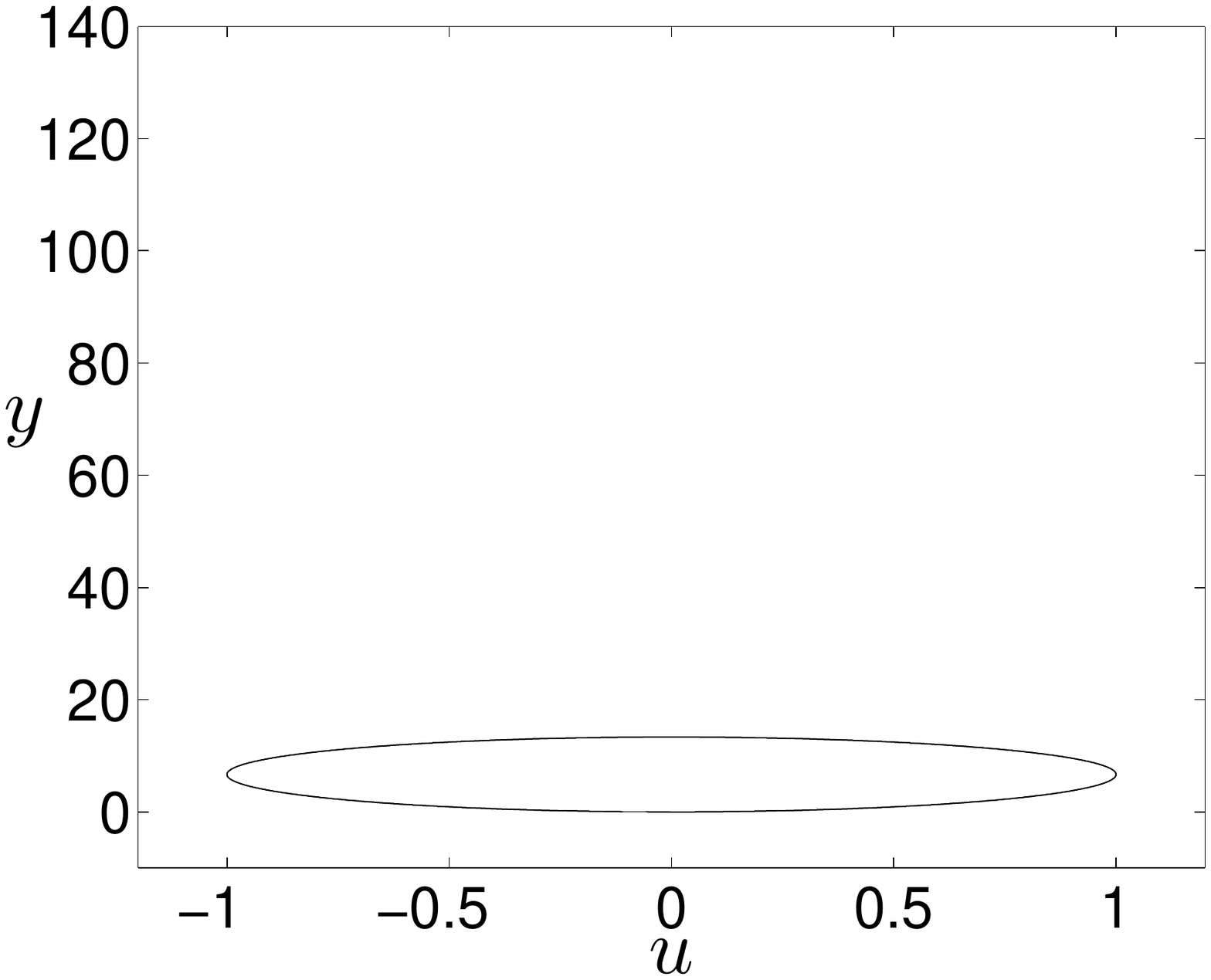}} 
  \subfloat[$\omega=0.001$]{\label{}\includegraphics[trim = 0mm 60mm 0mm 60mm, clip,width=0.26\textwidth]{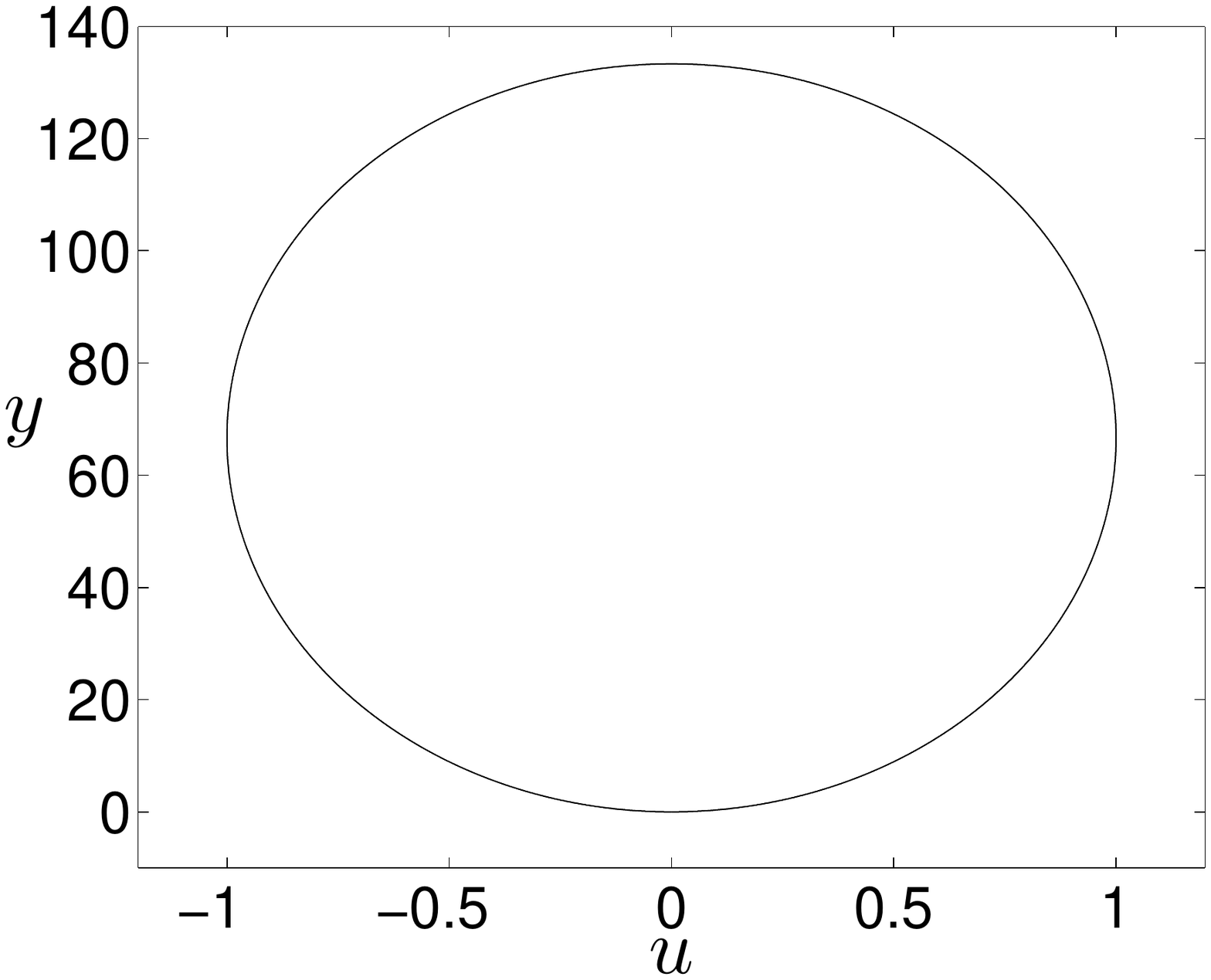}}
  \caption{\label{figlinearspringtrolleyloopmultieq}\small Input-output curves for the linear system in (\ref{eqspringsystemkis0}).  There are an infinite number of equilibria, $(a, 0)$. The initial condition is $(y,\dot y)=(0,0)$ and the input is $u(t)=\sin(\omega t)$ for various $\omega$.}
\end{figure}

 For arbitrary initial conditions $(y(0),\dot y(0))=(y_0,y_1)$,  the solution to (\ref{eqspringsystemkis0}) with $u(t)=\sin(\omega t)$ is 
\begin{align*}
y(t) 
&=y_0+\frac{y_1}{c}\left(1-e^{-ct}\right)-\frac{\sin(\omega t)}{\omega^2+c^2}\\
&+\frac{\omega^2+c^2-c^2\cos(\omega t)-\omega^2 e^{-ct}}{\omega c \left(\omega^2+c^2\right)}.
\end{align*}
It follows that 
\[
\lim_{\omega \rightarrow 0} y(t)=y_0+\frac{y_1}{c}(1-e^{-ct})
\]
and hence 
\[
\lim_{\underset{\mathlarger{\omega\rightarrow 0}}{t\rightarrow \infty}}y(t)=y_0+\frac{y_1}{c}.
\]
The steady-state  of (\ref{eqspringsystemfirstorderMultieq}) as $\omega \rightarrow 0$ depends on the initial conditions. Input-output curves are displayed in Figure~\ref{figlinearspringtrolleyloopmultieq} with $c=15$,  zero initial condition and $u(t)=\sin(\omega t)$ for various values of the frequency $\omega$. The figures indicate that a loop persists as $\omega$ approaches $0$. 

Note that the loop is a smooth circle, and does not display the same rapid transition as in Figures 3d and 4. The continuum of equilibrium points, rather than disconnected equilibrium points, may be the reason for the absence of sharp jumps. We will see in the next section that hysteresis in the Landau-Lifshitz equation behaves similarly.

\section{LANDAU-LIFSHITZ EQUATION}\label{secLLeq}
The Landau-Lifshitz equation is a coupled set of three nonlinear partial differential equations, 
\begin{subequations}\label{eqLLmodel}
\begin{align}
\dot{\mathbf m} &= \mathbf m \times  \mathbf m_{xx}-\nu\mathbf m\times\left(\mathbf m\times\mathbf m_{xx}\right)\label{eqLLGuoDing}\\
    \mathbf m_x(0,t)&=  \mathbf m_x(L,t)=\mathbf 0.\label{eqboundarycondition}\\
      \mathbf m(x,0)&=\mathbf m_0(x)\label{eqinitialcond}
\end{align}
where $\times$ indicates the vector cross product, $L>0$ is  length and $\nu >0$ a material parameter.  The solution to (\ref{eqLLmodel}) is 
\[
\mathbf m(x,t) = (m_1(x,t), m_2(x,t),m_3(x,t))
\]
for $x\in[0,L]$and 
\begin{align*}
\mathbf m \times  \mathbf m_{xx} =&\left(\right.m_2m_3^{\prime\prime} -m_3m_2^{\prime\prime}, \\
& \left.-m_1m_3^{\prime\prime} +m_3m_1^{\prime\prime}, m_1m_2^{\prime\prime} -m_2m_1^{\prime\prime}\right).
\end{align*}
The prime notation indicates differentiation with respect to $x$. Solutions to (\ref{eqLLmodel}) are defined on $\mathcal L_2^3 = \mathcal L_2 ([0,L]; \R^3).$   Existence and uniqueness of solutions to (\ref{eqLLmodel}) is shown in  \cite{Carbou2001}.    Equation~(\ref{eqLLmodel}) has been investigated in \cite{Alouges1992, Carbou2001}, \cite[Section~6.3.1]{Guo2008}.  

On a molecular level, ferromagnets are divided into regions called {magnetic domains}.  Each domain in a ferromagnet is magnetized to the same saturation, $M_s$. Mathematically, this means  $
  || \mathbf m(x,t)||_{2} =M_s$ where $||\cdot||_{2}$ is the Euclidean norm.  
  In much of the literature, $M_s$ is set to $1$ and the same convention is followed here; that is, 
\begin{equation}\label{eqconstraint}
  || \mathbf m(x,t)||_{2} =1
\end{equation} 
\end{subequations}
and the initial condition $||\mathbf m_0(x)||_2=1$.   The nonnegative parameter $\nu$ in (\ref{eqLLGuoDing}) is the damping parameter. It depends on the magnetization saturation and hence on the type of ferromagnet. 

\begin{theorem} \label{thmgeneraleqsol}  The set of  equilibrium solutions to the Landau-Lifshitz equation (\ref{eqLLGuoDing}) with Neumann boundary conditions  (\ref{eqboundarycondition}) is ${\mathbf m}(x)=\mathbf a$ where $\mathbf a \in \R^3$,  $||\mathbf a||_{2}=1. $
\end{theorem}

\noindent
{\em Proof:} It is shown in  \cite[Theorem~6.1.1]{Guo2008} that the equilibrium solutions to (\ref{eqLLGuoDing}) are
\[
{\mathbf m}(x)=\mathbf a \cos(kx)+\mathbf b \sin(kx)
\]
 where $k\in\mathbb R$ is a constant, $\mathbf a =(a_1,a_2,a_3)$ and $\mathbf b=(b_1,b_2,b_3)$ for $a_i,b_i \in \mathbb R$, $i=1,2,3$ constants with $||\mathbf a||_{2}=1,$ $||\mathbf b||_{2}=1$ and $\mathbf a^{\mathrm T} \mathbf b=0$.
Applying the boundary conditions in (\ref{eqboundarycondition}) implies $k=0$. 
  \hfill $\Box$
 
 Note that, as for the second-order linear example in (\ref{eqspringsystemfirstorderMultieq}), the Landau-Lifshitz equation has an infinite number of equilibria.  The set of equilibrium points of (\ref{eqLLmodel}) is
\begin{align}\label{equilibriumset}
E=&\{\mathbf a=(a_1,a_2,a_3) \in \mathbb{R}^3:\nonumber \\
&a_1,a_2,a_3 \mbox{ constants and }\mathbf a^\mathrm{T}\mathbf a=1 \}.
\end{align}

\begin{theorem}\label{thmlyapunov}  \cite[Proposition~6.2.1]{Guo2008}
Each point in the equilibrium set  $E$ (\ref{equilibriumset}) is stable in the $\mathcal L_2^3$-norm.
\end{theorem}

If we introduce the notion of an equilibrium set, then $E$ can be shown to be asymptotically stable \cite{Amenda_thesis}. 
 
Since the Landau-Lifshitz equation has an infinite number of stable equilibria, Definition~\ref{defmultiequilibrium} suggests it exhibits hysteresis. 
 
 Consider the Landau-Lifshitz equation with input $\mathbf u(t)$:
\[
\dot{ \mathbf m} = \mathbf m \times  \mathbf m_{xx}-\nu\mathbf m\times\left(\mathbf m\times\mathbf m_{xx}\right) +{\mathbf u}(t).
\]
 The input is chosen to be ${\mathbf u}(t)=\left(0.001\cos(\omega t),0,0\right)$; that is, there is a time-varying magnetic field in the $m_1$ direction and no magnetic field in the $m_2$ and $m_3$. 

\begin{figure}[h]
\vspace*{-0.3cm}
\centering\hspace*{-0.25cm}
    \subfloat[$\omega=1$]{ \includegraphics[trim = 0mm 60mm 0mm 60mm, clip,width=0.26\textwidth]{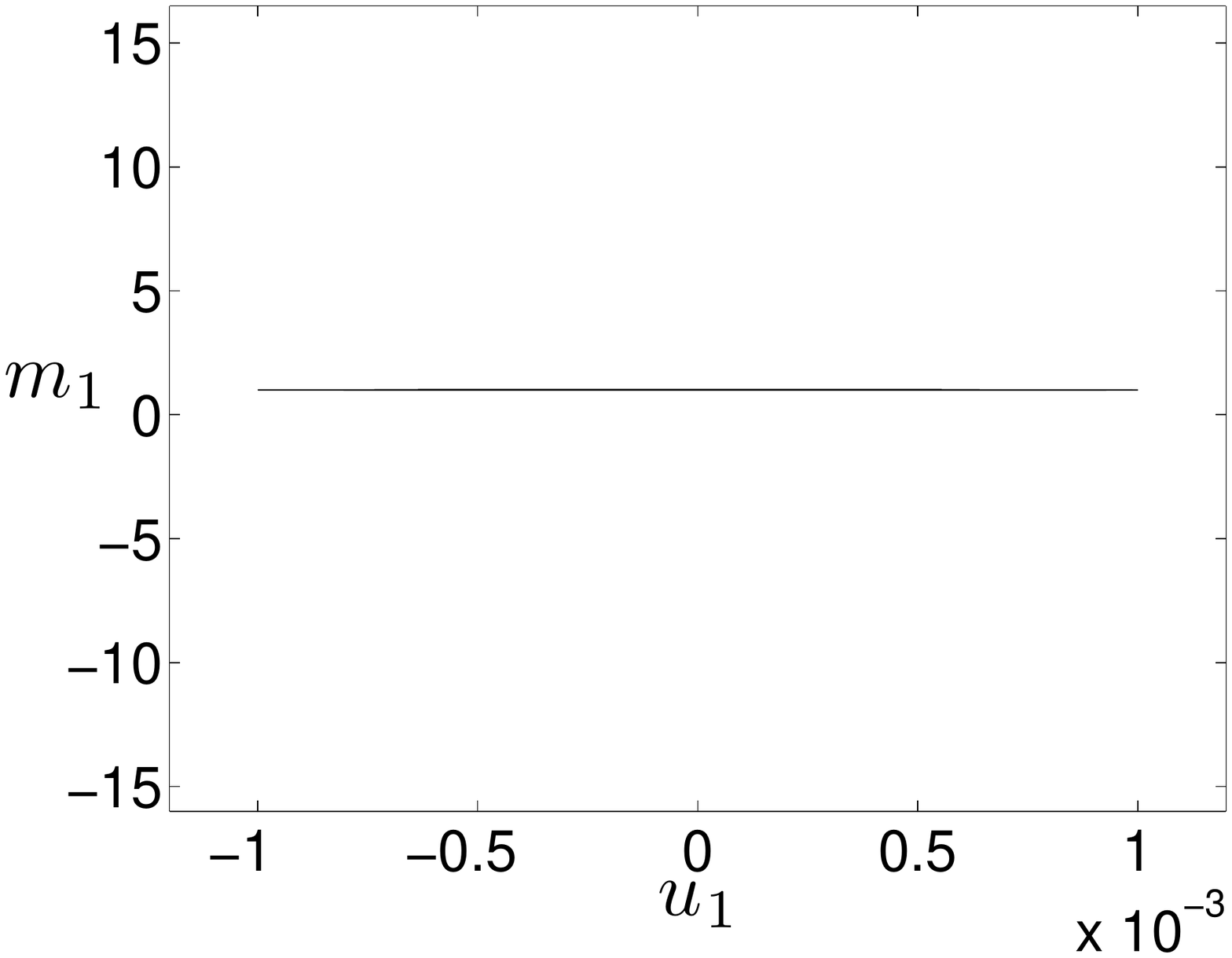}} 
  \subfloat[$\omega=0.1$]{\includegraphics[trim = 0mm 60mm 0mm 60mm, clip,width=0.26\textwidth]{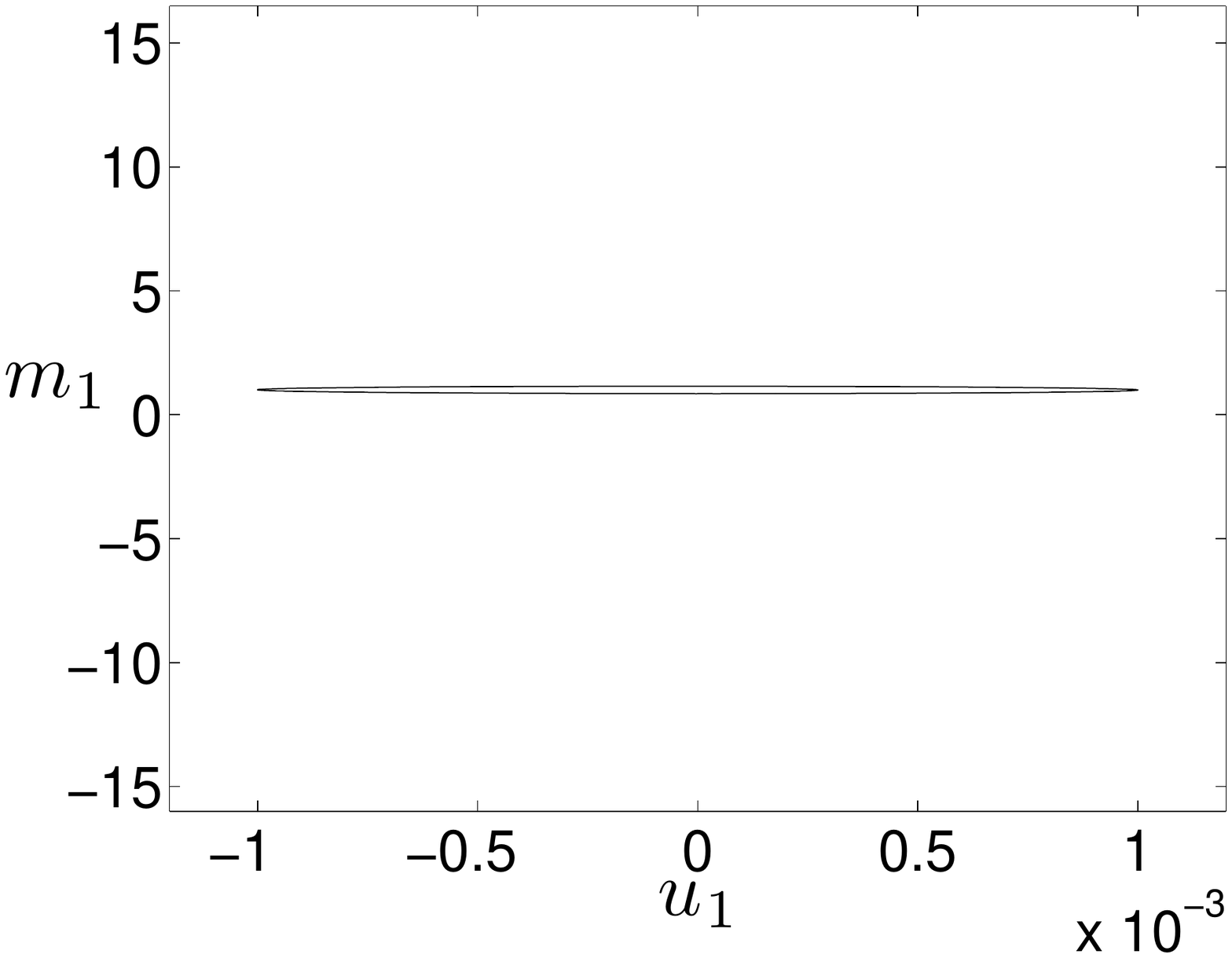}}\\
  \hspace*{-0.25cm}
  \subfloat[$\omega=0.01$]{ \includegraphics[trim = 0mm 60mm 0mm 60mm, clip,width=0.26\textwidth]{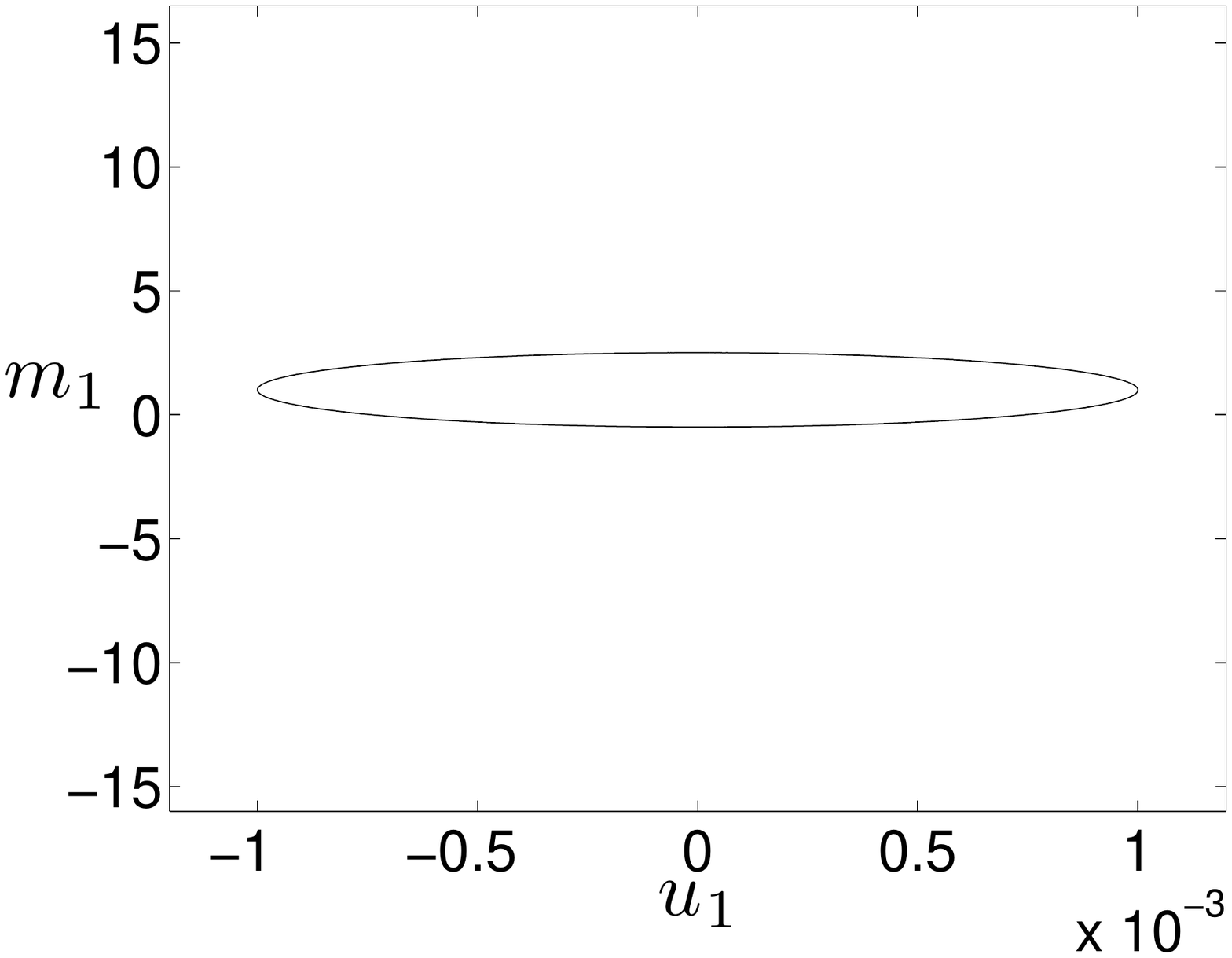}} 
    \subfloat[$\omega=0.001$]{ \includegraphics[trim = 0mm 60mm 0mm 60mm, clip,width=0.26\textwidth]{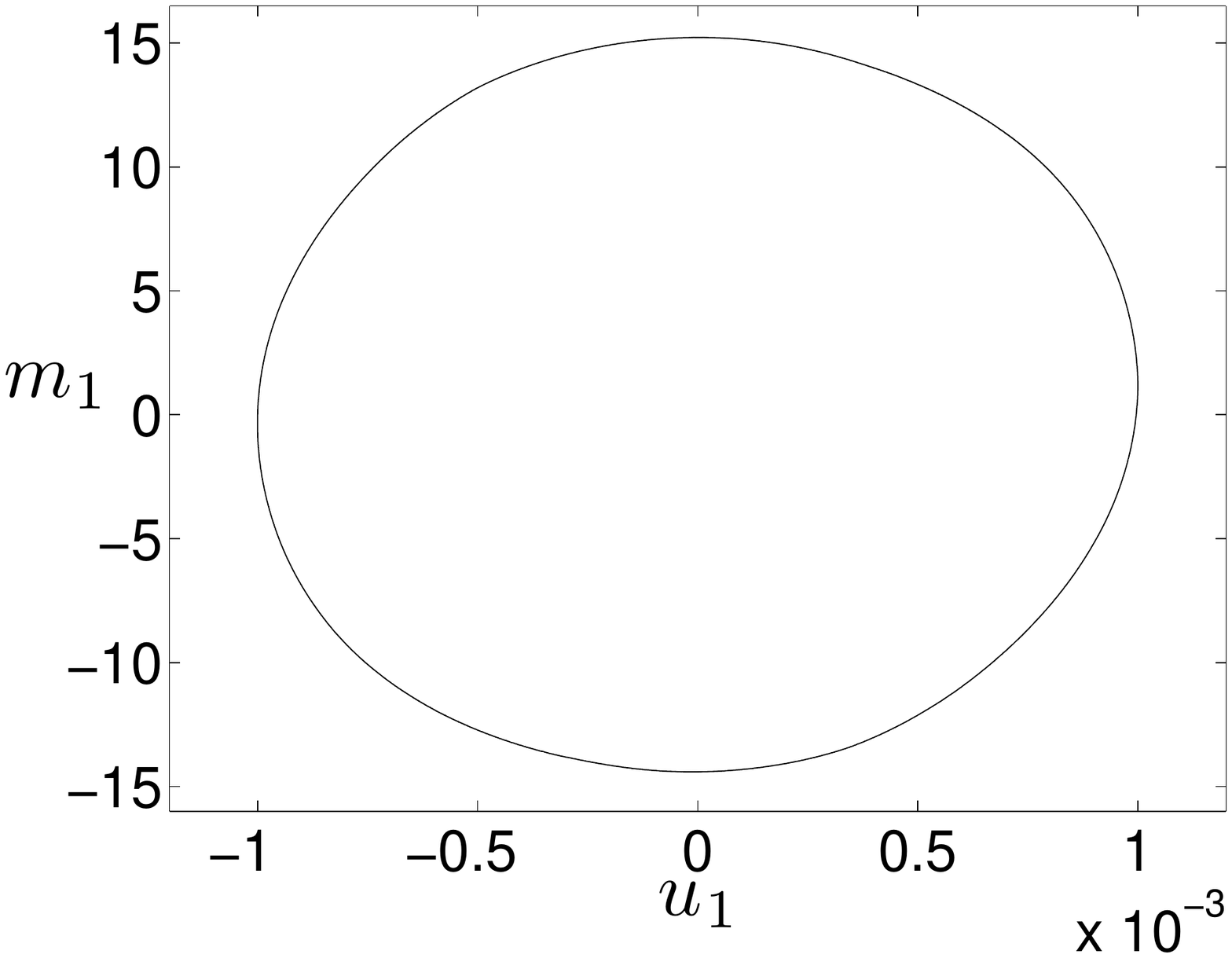}}    \vspace*{-0cm}
\caption{\label{figLoopFrequencyVaried} \small Hysteresis loops for $m_1(x,t)$ of the (nonlinear) Landau-Lifshitz equation with $x=0.6$ and $\nu=0.02$. The input is $\mathbf u(t)=\left(0.001\cos(\omega t),0,0\right)$. The initial condition is $\mathbf m_0(x)=\left(1,0,0\right)$.}
\end{figure}

The hysteresis loops for $m_1(x,t)$ with $x$ fixed is illustrated in Figure~\ref{figLoopFrequencyVaried}.  The initial condition is $\mathbf m_0(x)=\left(1,0,0\right)$, $\nu=0.02$ and $L=1$.   It is clear from Figure~\ref{figLoopFrequencyVaried} that the input-output curves  exhibit persistent loops as the frequency of the input approaches zero. In a similar manner, $m_2$ and $m_3$ also exhibit persistent loops.  Thus, these simulations indicate that  the Landau-Lifshitz equation satisfies Definition~\ref{defBernstein}.

Note that the hysteresis loops in Figure \ref{figLoopFrequencyVaried}  have the same appearance as those in Figure \ref{figlinearspringtrolleyloopmultieq} for the second-order linear example (\ref{eqspringsystemfirstorderMultieq}). In both examples there is a continuum of  equilibrium points.

To obtain the linear Landau-Lifshitz equation, equation (\ref{eqLLGuoDing}) is rewritten in semilinear form,
\begin{equation}\label{eqLLGuoDingTwo}
\frac{\partial \mathbf m}{\partial t} =\nu  \mathbf m_{xx}
 +\mathbf m \times  \mathbf m_{xx}
+\nu\left| \left|  \mathbf m_{x}\right| \right|_{2}^2 \mathbf m,
\end{equation}
using equation~(\ref{eqconstraint}) and properties of cross products.  Then substitute $\mathbf m (x,t) = {\mathbf a}+ \mathbf z (x,t)$
 where $\mathbf a \in E$ is an equilibrium of (\ref{eqLLmodel}) and $\mathbf z \in \mathcal L_2^3$ is a small perturbation into (\ref{eqLLGuoDingTwo}).
Defining
$$A\mathbf z=  \nu\mathbf z_{xx}+\mathbf a \times \mathbf {z}_{xx}$$
with domain
\begin{align*}
D(A)= \{& \mathbf z: \mathbf z\in \mathcal L_2^3, \,\,  \mathbf  z_{x} \in \mathcal L_2^3, \, \,  \mathbf z_{xx}^{} \in \mathcal L_2^3,\\
& \mathbf  z_{x}(0,t)=\mathbf 0=\mathbf  z_{x}(L,t)\},
\end{align*}
the Landau-Lifshitz equation linearized about an equilibrium $\mathbf a$ is
 \begin{subequations}\label{eqstatespaceform}
\begin{align}\dot{ \mathbf z}&=A\mathbf z,\quad \mathbf z(0)=\mathbf z_0 \\
\mathbf  z_{x}(0,t)&=\mathbf 0=\mathbf  z_{x}(L,t) .
\end{align}
\end{subequations} 

The operator $A$ can be shown to generate an analytic semigroup \cite{Amenda_thesis}. Therefore, (\ref{eqstatespaceform}) is a well-posed system. 

\begin{theorem}\label{thmlinLLstable}
Any constant $\mathbf c \in \R^3$ is a stable equilibrium of (\ref{eqstatespaceform}).
\end{theorem}

\noindent
{\em Proof:}
A brief outline is given here; details can be found in \cite{Amenda_thesis}.
Since $A$ generates an analytic semigroup, the spectrum determined growth assumption is satisfied and so the eigenvalues of $A$ determine the stability of the linear system (\ref{eqstatespaceform}) \cite[Section~5.1]{Curtain1995}, \cite[Section~3.2]{Luo1999}.

It is clear that any constant function $\mathbf c$ is an equilibrium of (\ref{eqstatespaceform}). Let $\lambda \in \mathbb C$. The eigenvalue problem of (\ref{eqstatespaceform}) is $\lambda\mathbf v=A\mathbf v$ and  
boundary conditions $\mathbf v_x(0)=\mathbf v_x(L)=\mathbf 0$.
where $\mathbf v \in \mathcal L_2^3$. 
Solving, the eigenvalues of (\ref{eqstatespaceform}) are $\lambda_1=0$ and for $n\in \mathbb Z$,
  \begin{align*}
 \lambda_2^{+,-} &= \frac{-(1+2n)^2\pi^2\nu}{L^2}\pm i\frac{(1+2n)^2\pi^2}{L^2}, \\
         \lambda_3 &=  \frac{-(1+2n)^2\pi^2\nu}{L^2}, \\
     \lambda_4^{+,-} &= \frac{-(2n)^2\pi^2\nu}{L^2}\pm i\frac{(2n)^2\pi^2}{L^2},\\
            \lambda_5 &= \frac{-(2n)^2\pi^2\nu}{L^2}.
               \end{align*}
Since all the eigenvalues have nonpositive real part, the equilibria of (\ref{eqstatespaceform}) are stable. 
\hfill $\Box$
 
\begin{figure}[h]\vspace*{-0.3cm}
\centering 
\hspace*{-0.25cm}
    \subfloat[$\omega=1$]{ \includegraphics[trim = 0mm 60mm 0mm 60mm, clip,width=0.26\textwidth]{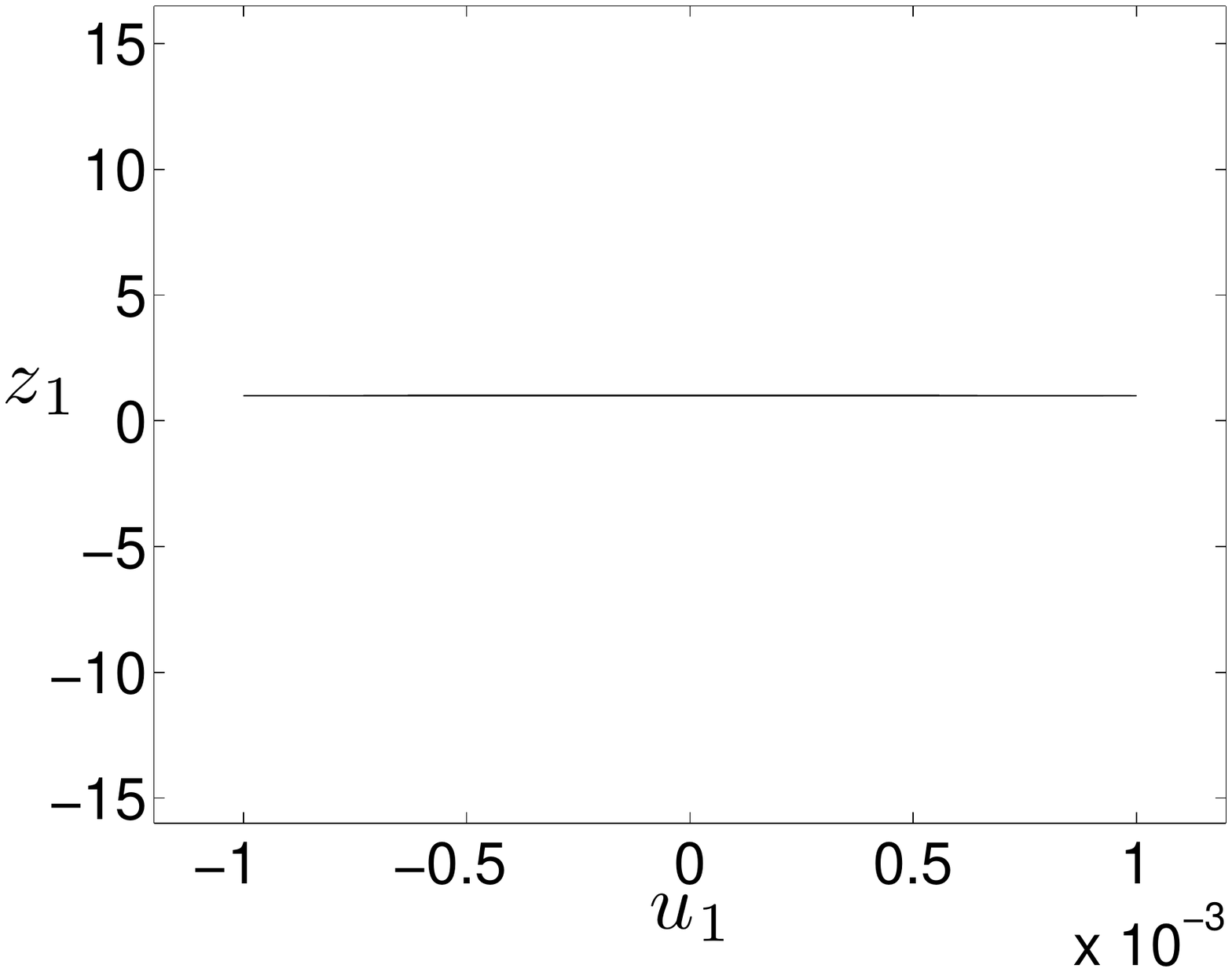}} 
  \subfloat[$\omega=0.1$]{\includegraphics[trim = 0mm 60mm 0mm 60mm, clip,width=0.26\textwidth]{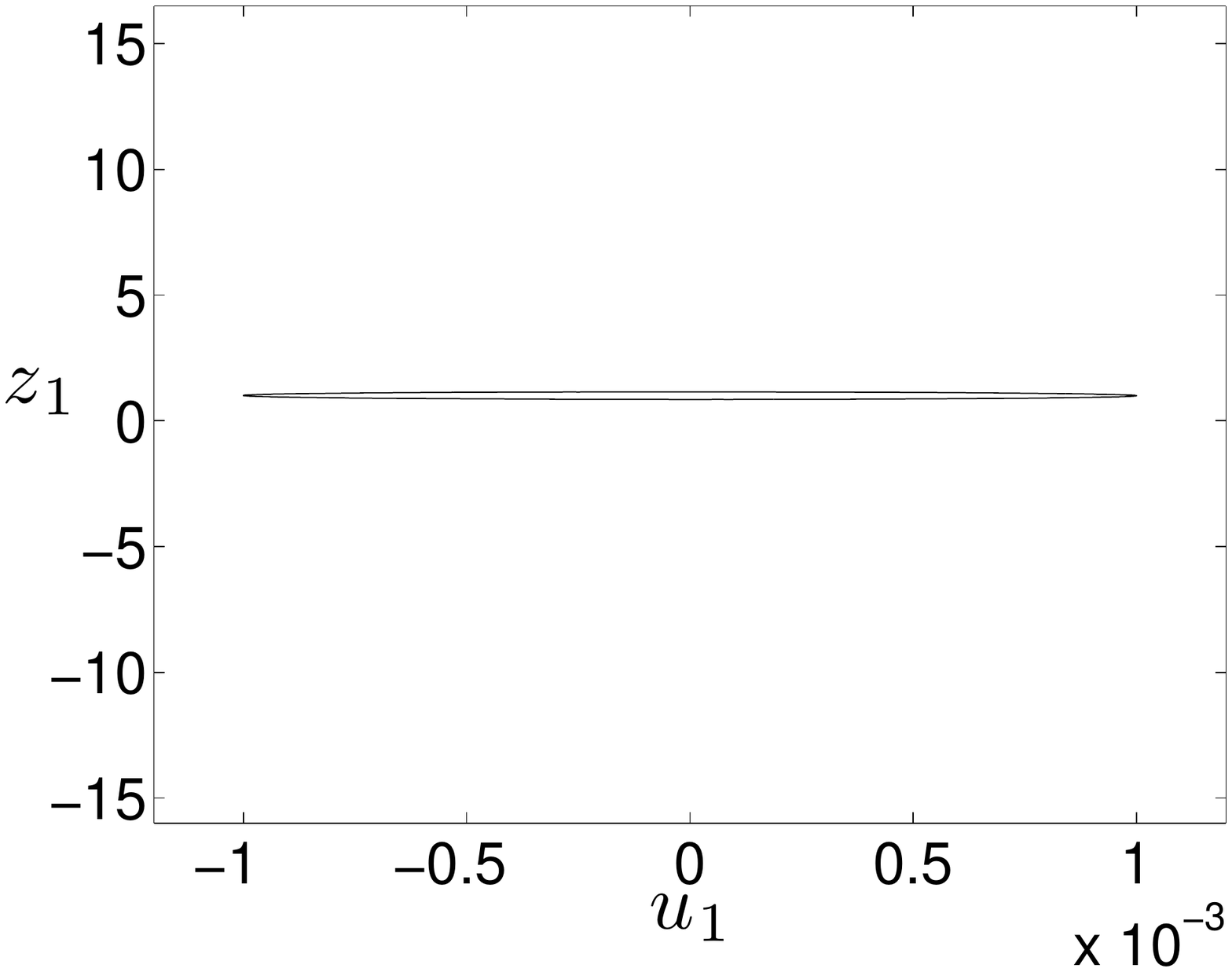}}\\
  \hspace*{-0.25cm}
  \subfloat[$\omega=0.01$]{ \includegraphics[trim = 0mm 60mm 0mm 60mm, clip,width=0.26\textwidth]{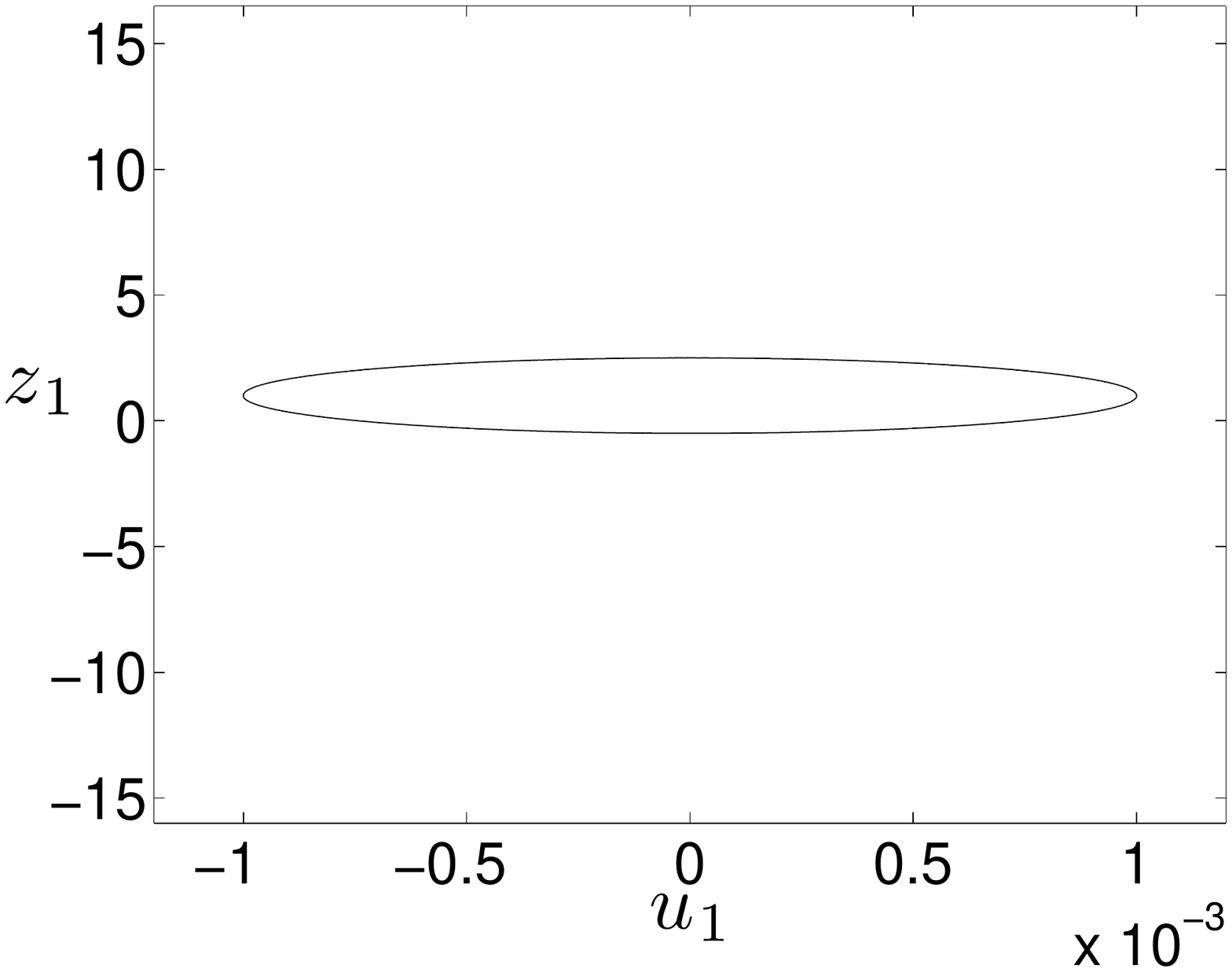}} 
    \subfloat[$\omega=0.001$]{ \includegraphics[trim = 0mm 60mm 0mm 60mm, clip,width=0.26\textwidth]{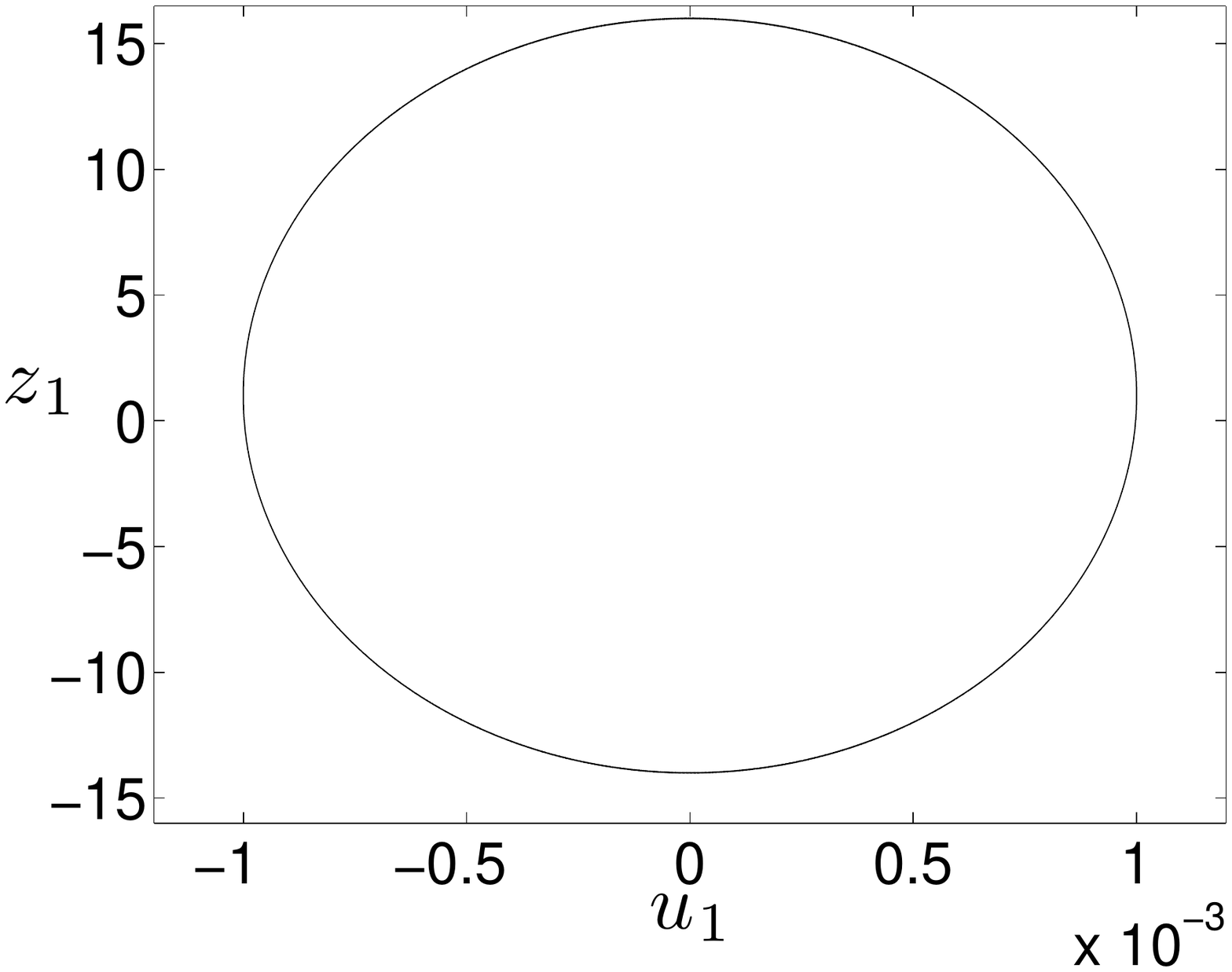}} 
\caption{\label{figlinearhysteresisloopM1} \small Hysteresis loops for $z_1(x,t)$ of the linear Landau-Lifshitz equation with $x=0.6$ and $\nu=0.02$.  The linearization is at $\mathbf a=(1,0,0)$. The input is $\mathbf u(t)=\left(0.001\cos(\omega t),0,0\right)$ and the initial condition is $\mathbf z_0(x)=\left( 1, 0,0 \right)$.}
\end{figure}

Based on the result in Theorem~\ref{thmlinLLstable},  Definition~\ref{defmultiequilibrium}  indicates that the linear Landau-Lifshitz equation exhibits hysteresis. We now look at  periodic inputs to determine whether persistent loops exists and Definition~\ref{defBernstein} is satisfied.  Again, $\nu=0.02$,  $L=1$,  and the same periodic input is applied as for the nonlinear Landau-Lifshitz equation.  Figure~\ref{figlinearhysteresisloopM1} shows the input-output curves for the linear Landau-Lifshitz equation with $\mathbf a=(1,0,0)$ and initial condition  $\mathbf z_0(x)=\left( 1, 0,0 \right)$.   A loop does persist  as the frequency of the input approaches zero. Similar plots are obtained when control is on the second and third components. The hysteresis loops in Figure~\ref{figlinearhysteresisloopM1}  are similar in shape to  the nonlinear Landau-Lifshitz equation depicted in Figure~\ref{figLoopFrequencyVaried} and the example in~(\ref{eqspringsystemfirstorderMultieq}), each of which have a continuum of equilibria.

\section{CONCLUSIONS}\label{secDiscussions}

The results here indicate that existence of multiple stable equilibria is crucial for systems to be hysteretic.   The linear and nonlinear Landau-Lifshitz equations  exhibit hysteretic behaviour. Hysteresis was also shown in  a second-order linear differential equation. Since hysteretic behaviour was displayed by several linear systems, this suggests that multiple stable equilibria, not nonlinearity, are the key factors leading to the display of hysteresis in systems. 
As well,  an absence of jumps in the hysteresis loops were observed. Further investigation is required to determine whether this is always the case for systems with a zero eigenvalue. 

 In \cite{Chow2014}, the authors discuss controlling hysteresis in the Landau-Lifshitz equation; that is controlling between equilibrium points.  





\bibliographystyle{is-abbrv}
\bibliography{ref}

\end{document}